\documentclass{article}
\usepackage{setspace}
\usepackage{amsmath,amsthm,verbatim,amssymb,amsfonts,amscd, graphicx}
\usepackage{graphics}
\providecommand{\keywords}[1]{\textbf{\textbf{Keywords:}} #1}
\usepackage{lipsum}

\usepackage{enumerate}
\usepackage{color}
\usepackage{fix-cm}
\usepackage{mathrsfs}
\usepackage[paper=a4paper, dvips, top=3cm, left=2.7cm, right=2.7cm, foot=1.3cm, bottom=3.2cm]{geometry}
\newtheorem{thm}{Theorem}[section]
\newtheorem{lemma}[thm]{Lemma}
\newtheorem{cor}[thm]{Corollary}
\newtheorem{prop}[thm]{Proposition}

\begin{document}
\title{One-dimensional System Arising in Stochastic Gradient Descent}

\author{Kostis Karatapanis}

\date{}

\maketitle
\begin{abstract}
	We consider SDEs of the form $dX_t = |f(X_t)|/t^{\gamma} dt+1/t^{\gamma} dB_t$, where $f(x)$ behaves comparably to $|x|^k$ in a neighborhood of the origin, for $k\in [1,\infty)$. We show that there exists a threshold value $:=\tilde{\gamma}$ for $\gamma$, depending on $k$, such that when $\gamma \in (1/2, \tilde{\gamma})$ then $\mathbb{P}(X_n\rightarrow 0) = 0$, and for the rest of the permissible values $\mathbb{P}(X_n\rightarrow 0)>0$. The previous results extend for discrete processes that satisfy $X_{n+1}-X_n = f(X_n)/n^\gamma +Y_n/n^\gamma$. Here, $Y_{n+1}$ are martingale differences that are a.s. bounded.
	
	This result shows that for a function $F$, whose second derivative at degenerate saddle points is of polynomial order, it is always possible to escape saddle points via the iteration $X_{n+1}-X_n =F'(X_n)/n^\gamma +Y_n/n^\gamma$ for a suitable choice of $\gamma$.
\end{abstract}
\bigskip
	\begin{section}{Introduction}\label{Intro}
			Let $F:\mathbb{R}^d \rightarrow \mathbb{R}^d,~d \geq 1$ be a vector field. For much of what follows $F$ arises the as gradient of a potential function $V$, namely $V:\mathbb{R}^d \rightarrow \mathbb{R}$, and $F=-\nabla V$. Now, we define a system driven by
			\begin{equation}\label{eq:King}
			X_{n+1}=X_n+ a_n(F(X_n) +\xi_{n+1})).
			\end{equation}
			To elaborate on the parameters, let $\mathscr{F}_n$ be a filtration, then $a_n,\xi_n$ are adapted, and $\xi_n$ constitute martingale differences i.e. $E(\xi_{n+1}|F_n) = 0$. For the purposes of this introduction we will simplify and assume, without any great loss of abstraction, that $a_n$ is deterministic and is either a constant number or it is converging to zero comparably to $n^{-\gamma},~\gamma \in (1/2,1]$. Also, some additional assumptions on the noise are usually required: one is a boundedness restraint, that is we assume the existence of a constant $M$, such that $|\xi_n|\leq M$ a.s.; and secondly we want $\xi_n$ to be isotropic i.e. $\mathbb{P}( (\theta \cdot \xi_n)^{+}>\delta    )>\delta$ for any unit direction $\theta \in \mathbb{R}^d$.
			This versatile system is well studied, and it arises naturally in many different areas. In machine learning and statistics $\eqref{eq:King}$ can be used as a powerful tool used for quick optimization and statistical inference, among other uses. Furthermore, many urn models are represented by \eqref{eq:King}. These processes play a central role in probability theory due to their wide applicability in physics, biology and social sciences; for an comprehensive exposition on the subject see~\cite{pm2007}.
			
			In machine learning, processes satisfying \eqref{eq:King} appear in stochastic gradient descent (SGD). First, to provide context, let us briefly introduce the gradient descent method (GD) and then see why SGD
			arises naturally from there. The GD is an optimization technique which finds local minima for a potential function $V$ via the iteration
			\begin{equation}
			x_{n+1}-x_{n}=-\eta_n \nabla V(x_n), 
			\end{equation}
			in many applications we take $\eta_n$ to be a positive and constant.
			The previous method when applied to non-convex functions has the shortcoming that it may get stuck
			near saddle points, i.e. points where the gradient vanishes, that are neither local minima nor local
			maxima, or locate local minimuma instead of global ones. The former issue can be resolved by adding noise into the system, which, consequently, helps in pushing the particle downhill and eventually escaping saddle points. For the latter, in general, avoiding local minima is a difficult problem (\cite{ConvergenceToGM1} and \cite{ConvergenceToGM2}), however, fortunately, in many instances finding local minima is satisfactory. Recently, there have been several problems of interest where this is indeed the case, either because all local minima are global minima (\cite{Ge15} and \cite{Dict}), or in other cases local minima provide equally good results as global minima~\cite{Chromo}. Furthermore, in certain applications saddle points lead to highly sup-optimal results (\cite{NonOptimalSaddle1} and \cite{NonOptimal2}).
			
		As described in the previous paragraph escaping saddle points when performing SGD is an important problem. In the literature the saddle problem, when non-degeneracy conditions are imposed, is well understood. Results showing that asymptotically SGD will escape saddle points, date back to works of Pemantle~\cite{MR1055778}, and, more recently,~\cite{Lee16} where they prove random initialization guarantees almost sure convergence to minimizers. Having established asymptotic convergence, subsequently, led to results on how this can be done efficiently~\cite{Lee16}.
			
			Processes satisfying \eqref{eq:King}, when $a_n$ goes to zero, are known as stochastic approximations after~\cite{MR0042668}. These processes have been extensively studied since then~\cite{kushner2003stochastic}. An important feature is that the step size $a_n$, satisfies  $$\sum_{n\geq 1} a_n =\infty \text{~and~}  \sum_{n\geq 1} a_n^2 <\infty.$$ This property balances the effects of the noise in the system, so that there is an implicit averaging that, eventually, eliminates the effects of the noise. The previously described system hence behaves similarly to the mean flow: the ODE whose right-hand side corresponds to the expectation of the driving term. The previous heuristic, can helps us identify the support $S$ of the limiting process $X_\infty:=\lim_{n \rightarrow \infty} X_n$, in terms of the topological properties of the dynamical system $\frac{\text{d }\!\!X_t   }{ \text{d }\!\! t}=F(X_t) $. More specifically, in most instances, one can argue that attractors or ``strict" saddle are in $S$, whereas repellers are not (see \cite{kushner2003stochastic}). However, there has not been a systematic approach finding when a degenerate saddle point, i.e. a point that is neither an attractor nor a repeller, belongs in $S$.
			
			Stochastic approximations arise naturally in many different contexts. Some early results were published by ~\cite{Ruupert}, and ~\cite{Polyak}. There, they dealt with averaged stochastic gradient descent (ASGD) arising from a strongly convex potential $V$ with step size $n^{-\gamma},~\gamma \in (1/2,1]$. In their work they proved that one can build, with proper scaling, consistent estimators $\tilde{x}_n$ (for the $\arg \min(V)$) whose limiting distribution is Gaussian.  In learning problems, a modified version of ASGD ~\cite{Sasha} provides convergence rates to global minima of order $n^{-1}$. Additionally, many classical urn processes can be described via \eqref{eq:King}, where $a_{n}$ is of the order of $n^{-1}$. Certain effort is being placed in understanding the support of limiting process $X_{\infty}$. In specific instances, the underlying problem boils down to understanding an SGD problem: characterizing the support of $X_\infty$ in terms of  the class of critical points of the corresponding potential $V$. For a comprehensive exposition on urn processes see~\cite{pm2007}.

			From the previous discussion some fundamental questions of interest regarding \eqref{eq:King} are: 
			\begin{enumerate}
				\item Does  $X_n$ converge?
				\item When does $X_n$ converge to minima (local), consequently avoiding saddle points?
				\item When does $X_n$ converge to global minima?
				\item How fast  does $X_n$ converge to local minima?
			\end{enumerate}
			When $F$ arises from a potential function $V$, the first question is for the most part settled: the limit of the process converges, and it is supported on a subset of the set of critical points of $V$.
			
			Here, our primary focus will be understanding the second question in a one-dimensional setting. In the literature there are many results of this type. However, as already mentioned, the vast majority of them require the saddle points to satisfy certain non-degeneracy conditions. In fact, non-degenerate saddle point will never be in the support of $X_\infty$. Interestingly enough, the previous conclusion is not always valid for degenerate ones~\cite{MR1055778}. However, we show that for any $V$, under some mild conditions, we can find $\gamma$ such that saddle points do not belong in $S$. Hence, we demonstrate that implementing SGD, by adding enough noise, gives the desired asymptotic behavior even in the degenerate case. 
						
			Let us introduce some motivating background originating from urn model theory. In paper~\cite{TOUCH1980}, they consider a random variable $X_n$ taking values in $(0,1)$, which we interpret as counting the percentage of the red balls out of $n$ balls in accordance to Polya's urn model. Recursively define $X_{n+1}$ to be $\dfrac{nX_n +1}{n+1}  $ with  probability $f(X_n)$, and $X_{n+1} = \dfrac{nX_n}{n+1}  $ with probability $1-f(X_n)$. The main result of~\cite{TOUCH1980} is that $X_n$ converges to a random variable $X$, whose range is a subset of $C= \{p| f(p)=p    \}$; moreover for all points $p$ such that $f'(p)<1$ or ($f'(p)>1$), we have $\mathbb{P} (X=p     )>0$ ($\mathbb{P}(X=p)=0$).
			
			This process fits the general form \eqref{eq:King}. Indeed, we may rewrite $X_n$ in the following form
			$ X_{n+1}- X_n= A_n+ Y_n$ where $Y_n$ is the martingale
			$$
			Y_n = \left \{
			\begin{array}{ccc}
			\frac{1- f(X_n)}{ n+1} ,& \text{with\,probability\,}f(X_n)\\
			& &,\text{~and~} A_n= \frac{f(X_n) - X_n}{n+1}.\\ 
			\frac{- f(X_n)}{ n+1},& \text{with\,probability\,}1-f(X_n)
			\end{array} 
			\right .
			$$
			Define 
			$
			g_n = \left \{
			\begin{array}{cc}
			1- f(X_n) ,& \text{with\,probability\,}f(X_n)\\
			- f(X_n),& \text{with\,probability\,}1-f(X_n)
			\end{array} 
			\right . .
			$			
			
			Then, the SDE becomes 
			$   X_{n+1}- X_n=\frac{f(X_n) - X_n}{n+1} + \dfrac{g_n}{n+1}=\frac{f(X_n) - X_n}{n+1} + \dfrac{\Theta (1)}{n+1}     $, when $f(X_n)$ is bounded away from $\{0,1\}$. We have already mentioned that $X_n$ can only converge to points $p$, such that $f(p)=p,\, f'(p)<1$. The idea is that the condition $f'(p)<1$ implies that $f(X_n) - X_n $ is positive when $ X_n \in ( p-\delta ,p )$ and negative when $ X_n \in ( p ,p + \delta )$. Therefore, $A_n$ pushes $X_n$ towards $p$, when $ X_n$ lies in a neighborhood of $ p$, and since $|X_{n+1} - X_n| = O(1/n)$, the process $(X_n)_{n\geq 0}$ may eventually get trapped in the neighborhood. Consequentially, as $p$ is the sole point in the neighborhood that belongs in $C$, the convergence follows. For a detailed discussion of these results see chapter 2.4 in ~\cite{pm2007}.
		
			The previous analysis establishes that the support of the limiting process is exactly the set of fixed points of $f$ (critical points of the corresponding potential).More precisely, it will avoid local maxima with probability $1$, and it will converge to a local minimum with some positive probability. However, at that point in history, it was unknown whether $X_n$ can converge to saddle points. Later Pemantle, with his work~\cite{MR1055778}, settled this; giving explicit conditions, and surprisingly depending on the local behavior of $f$, the process  may or may not converge there. Next, we will define a quantity which we will need in the next paragraph. 
			Let $Z_{n,m} = \sum_{i=n}^{m-1} Y_i$, so
			$E (Z_{n,m}^2 ) \leq \sum _{i \geq n } \frac{1}{(i+1)^2} \sim \frac{1}{n}$. The last equation, after taking $m\rightarrow \infty$, is called the remaining variance for the process $X_n$, and it measures how much $X_n$ can potentially deviate from the ``mean flow" by the influence of future noise.
		
			We will give the high level intuition, in qualitative terms, utilizing objects already described, namely the mean flow and the remaining variance. It is clear that the occurrence of convergence or non convergence to a point $p$, depends on the behavior of the process $X_n$ when lying in the stable trajectory. Now, for simplicity, we assume the stable trajectory lies in a left neighborhood of $p$ namely $(p-\delta,p)$, and recalling that $p$ is a saddle point $(p,p+\delta)$ realizes the unstable trajectory. Consequently, assume $X_n$ is moving towards $p$. The notion of the expected rate of convergence $\mathfrak  {o}_1(n):=(X_n-p)$ is encapsulated by the mean flow which can be explicitly computed. To continue further, as promised, we need to introduce
			$\mathfrak  {o}_2(n)= \sqrt{E(Z_{n,\infty}^2   )}$ the order of the square root of the remaining variance. When $\mathfrak  {o}_1(n)=o( \mathfrak  {o}_2(n)  )$, in every instance where $X_n$ behaves as expected, with h.p. $X_n$ will be pushed, by the remaining noise, to the unstable trajectory i.e. $X_{n+k} \in (p,p+\delta)$ for some $k>0$. Whenever this happens $X_{n+k}$ may fail to return to $(p-\delta,p)$ with some positive fixed probability. Finally, by Borel-Cantelli the process will not converge to $p$ with probability $1$. Similarly, we can argue that when $\mathfrak  {o}_2(n)=o( \mathfrak  {o}_1(n)  )$, $X_n$ will converge to $p$ with some positive probability. To elaborate, the probability that $X_n$ will escape the stable trajectory is decaying rapidly whence by Borel-Cantelli, in the event that $X_n$ behaves as expected, the process will fail to visit the unstable trajectory, thereby establishing the convergence of $X_n$ to $p$. 
			\begin{subsection}{Results for the continuous model}
				We proceed by transitioning to a continuous model. For that purpose we need a potential, a step size, and a noise. However, it is natural to consider, without the need to contemplate, a process defined by
				\begin{equation}\label{eq:GenIntro}
				{ \rm d  }L_t = \frac{f(L_t)}{t^{\gamma}} { \rm d}t + \frac{1}{t^\gamma}   {\rm d} B_t,~\gamma \in (1/2,1]. 
				\end{equation}
				We assume that $f(0)=0$, and $f$, is otherwise positive in a neighborhood $\mathcal{N}$ of zero. What we wish to investigate is whether $L_t$ will not converge to $0$ with probability $1$, or if it will converge there with some positive probability. The answer to these questions depends only on the local behavior of $f$ on $\mathcal{N}$. 
			
				The main non-convergence result is the following:
				\begin{thm}\label{thm:gennonconv} Suppose that $\mathcal{N}$ is a neighborhood of zero.
					Let $(L_t)_{t\geq 1}$ be a solution of ~\eqref{eq:GenIntro}, where $f(x)$ is Lipschitz. We distinguish two cases depending on $f$ and the parameters of the system
					\begin{enumerate}
						\item \label{enum:nonconv1} $k|x| \leq f(x) $, $k> \frac{1}{2}$ and $ \gamma =1$ for all $x\in \mathcal{N}$.
						
						\item \label{enum:nonconv2}$ |x|^k  \leq f(x)\,$, $\frac{1}{2}+\frac{1}{2k} \geq \gamma$ and $k>1$ for all $x\in \mathcal{N}$.
					\end{enumerate}
					If either \ref{enum:nonconv1} or \ref{enum:nonconv2} hold, then $\mathbb{P}(L_t\rightarrow 0)=0$.
				\end{thm}
			     In the first part of the theorem, the result holds even in the case $k=\frac{1}{2}$, however the proof is omitted to avoid repetitiveness.  
			    In part \ref{enum:nonconv1}, we have only considered $\gamma=1$ since that is the only critical case, namely for $\gamma<1$ the effects of the noise would be overwhelming and for all $k$, we would obtain $\mathbb{P}(L_t\rightarrow 0)=0$.			    
			    
			    We now state the main convergence theorem:
			    \begin{thm}\label{thm:genconv} Suppose that $\mathcal{N}$ is a neighborhood of zero.
			    	Let $(L_t)_{t\geq 1}$ be a solution of~\eqref{eq:GenIntro}. We distinguish two cases depending on $f$ and the parameters of the system
			    	\begin{enumerate}
			    		\item \label{enum:conv1}$k_1|x|\leq f(x)\leq k_2|x|,~0<k_i<1/2$ and$~\gamma=1$ for all $x\in \mathcal{N}\cap(-\infty,0]$.
			    		
			    		\item \label{enum:conv2}$c|x|^k  \leq f(x)\leq|x|^k$,  $\frac{1}{2}+\frac{1}{2k} >\gamma~\text{and}~k>1$  for all $x\in \mathcal{N}\cap(-\infty,0]$.
			    		 
			    	\end{enumerate}
			    	If either \ref{enum:nonconv1} or \ref{enum:nonconv2} hold then $\mathbb{P}(X_t\rightarrow 0)>0$.
			    \end{thm}
				This is accomplished by first establishing the previous results for monomials i.e. $f(x)=|x|^k$ or $f(x)=k|x|$, which is done in sections \ref{sec:k=1}, and \ref{sec:k>1}. And we prove the stated theorems in section \ref{sec:ContGenCase}, by utilizing the comparison results found in section \ref{sec:comparisontheorems}.
			
				In section \ref{sec:k=1}, we deal with the linear case, i.e. when $f(x)=k|x|$. There, the SDE can be  explicitly solved, which simplifies matters to a great extend. Firstly, in subsection \ref{sbec:k=1Nconv}, we prove that when $k>1/2$, the corresponding process a.s. will not converge to $0$, which is accomplished by proving that it will converge to infinity a.s.. Secondly, in subsection \ref{sbec:k=1C} we show that process will converge to $0$ with some positive probability.
			
				In section \ref{sec:k>1}, we move on to the higher order monomials, i.e. $f(x)=|x|^k$. Here, we show that the process will behave as the ``mean flow" process $h(t)$ infinitely often, which is accomplished by studying the process $L_t/h(t)$. In subsection \ref{sbec:ContNConv}, the main theorem is that when $\frac{1}{2}+\frac{1}{2k}\geq\gamma$, then $L_t\rightarrow \infty$ a.s.. In the section \ref{sbec:contConv}, we show that when 
				$\frac{1}{2}+\frac{1}{2k}<\gamma$ holds, the process may converge to $0$ with positive probability.
				
				 Qualitatively, the previous constrains on the parameters are in accordance with our intuition. To be more specific, when $k$ increases, $f$ becomes steeper, which should indicate it is easier for the process to escape. When $\gamma$ decreases the remaining variance increases, hence we should expect that the process visits the unstable trajectory with greater ease, due to higher fluctuations. 
				
			\end{subsection}
		\begin{subsection}{Results for the discrete model}\label{sbec:DiscIntro}
				 The asymptotic behavior of the discrete processes is the expected one, depending on the parameters of the problem. Here, we study processes satisfying
			\begin{equation}\label{eq:genDiscNconv}
			X_{n+1} - X_n\geq\frac{f(X_n)}{n^{\gamma}} +  \frac{Y_{n+1}}{n^{\gamma}},\,\gamma \in (1/2,1),\,k\in (1,\infty), 
			\end{equation}
			or 
			\begin{equation}\label{eq:genDiscNconv2}
			X_{n+1} - X_n\leq\frac{f(X_n)}{n^{\gamma}} +  \frac{Y_{n+1}}{n^{\gamma}},\,\gamma \in (1/2,1),\,k\in (1,\infty), 
			\end{equation}
			where $Y_n$ are a.s. bounded, i.e. there is a constant $M$ such that $|Y_n|<M$ a.s., $E(Y_{n+1}|\mathscr{F}_n)=0$, and 
			$   E(Y_{n+1}^2|\mathscr{F}_n    )\geq l>0 $. 
			The main non-convergence theorem is the following
			\begin{thm}\label{thm:gennonconvdisc} Suppose that $\mathcal{N}$ is a neighborhood of zero.
				Let $(X_n)_{n\geq 1}$ solve~\eqref{eq:genDiscNconv}. If
				\begin{enumerate}				
					\item \label{enum:nonconv2disc}$ |x|^k  \leq f(x)\,$, $\frac{1}{2}+\frac{1}{2k} > \gamma$ and $k>1$ for all $x\in \mathcal{N}$
				\end{enumerate}			
				then $\mathbb{P}(X_n\rightarrow 0)=0$.
			\end{thm}
			For the convergence result the non-degeneracy condition $E(Y_{n+1}^2|\mathscr{F}_n)\geq l$ is replaced with the assumption stated in part \ref{assumption} of Theorem \ref{thm:genconvdisc2}.  
			\begin{thm}\label{thm:genconvdisc2} Let $\mathcal{N}=(-3\epsilon,3\epsilon)$ be a neighborhood of zero.
				Suppose $(X_n)_{n\geq 1}$ solve~\eqref{eq:genDiscNconv}. If
				\begin{enumerate}	
				\item	\label{assumption} There exist $-\epsilon _2>-3\epsilon$, $-\epsilon_1<-\epsilon $, such that for all $M>0$, there exists
					$n>M$ such that $\mathbb{P}(X_n \in (-\epsilon _2,-\epsilon_1    )            ) >0$.			
					\item \label{enum:conv2disc}$ 0<f(x)\leq |x|^k$, $\frac{1}{2}+\frac{1}{2k}< \gamma$ and $k>1$ for all $x\in \mathcal{N}$
				\end{enumerate}			
							Then $\mathbb{P}(X_n\rightarrow 0 )>0$
			\end{thm}
			The assumption imposed on $X_n$, part \ref{assumption} of Theorem \ref{thm:genconvdisc2}, says that the process should be able visit a neighborhood of the origin for large enough $n$. If this constraint is not imposed on the process, the previous result need not hold. For instance, the drift could dominate the noise, and, consequentially, the process may never reach a neighborhood of the origin with probability $1$. There are processes that naturally satisfy this property, such an example is the urn process defined in section \ref{Intro} (see \cite{MR1055778}).

			\end{subsection}
			\end{section}
			\begin{section}{Preliminary results}\label{sec:comparisontheorems}
				We will now prove two important lemmas, that we will be needed throughout. 
			Let $f:\mathbb{R}\rightarrow \mathbb{R}$, be Lipschitz such that for all $\epsilon>0$ there exists $c$ such that $f(x)>c,$ for all $x \in \mathbb{R}\setminus (-\epsilon,\epsilon)   $. Also, we define a continuous function $g:\mathbb{R}_{\geq 0}\rightarrow \mathbb{R}$, such that 
			$\int_{0}^{\infty} g^2(t) \text{d}t <\infty  $. Let $X_t$ that satisfies 
			\begin{equation}\label{eq:introinf}
			\text{d}X_t = f(X_t)\text{d}t +g(t)\text{d}B_t.
			\end{equation}
			\begin{lemma} \label{lem:ContLSup}
				$\limsup _{t \rightarrow \infty} X_t \geq 0  $ a.s..
			\end{lemma}
			\textbf{Proof:} We will argue by contradiction. Assume that $\limsup _{t \rightarrow \infty} X_t < 0   $, and pick $\delta >0   $ such that 
			$\limsup _{t \rightarrow \infty} X_t< -\delta$ with positive probability. Then there is a time $u$, such that $X_t \leq  -\delta$ for all $ t\geq u$. But this has as an immediate consequence that 
			$\int _{1}^{t} f(X_s) {\rm d}s \rightarrow \infty   $. However, since the process $G_t= \int _{1}^{t} g(s) {\rm d}B_s  $ has finite quadratic variation, i.e. $\sup_t \langle G_t \rangle=\int_{0}^{\infty} g^2(t) \text{d}t <\infty$,  $G_t$ stays a.s. finite. The last two observations imply that $X_t \rightarrow \infty$, which is a contradiction.
			$\hfill \blacksquare$

			\begin{lemma}\label{lem:ContLInf}
				$\liminf _{t \rightarrow \infty} X_t \geq 0  $ a.s..
			\end{lemma}
			\textbf{Proof:} We will again argue by contradiction. Assume that $\liminf _{t \rightarrow \infty} X_t < 0   $ on a set of positive probability. Take an enumeration of the pair of  positive rationals $(q_n,p_n)$ such that $q_n>p_n$. Now, define $A_{n}= \{     X_t \leq -q_n ~\text{i.o.}, X_t\geq  -p_n  ~\text{i.o.}   \}  $. Since $ \limsup _{t \rightarrow \infty} X_t \geq 0  $, we have $  \bigcup_{n\geq 0} A_{n}   =  \{ \liminf _{t \rightarrow \infty} X_t < 0           \}  $. Now, for $ t_1<t_2$ assume that $ X_{t_1} \geq -p_n    $ and 
			$ X_{t_2} \leq -q_n    $. Then, we see that $X_{t_2} -X_{t_1} \leq -q_n + p_n     $, however
			\begin{align*}
			X_{t_2} - X_{t_1} &= \int_{t_1}^{t_2} f(X_s) {\rm d} s  + \int _{t_1}^{t_2} g(s){\rm d}B_s   \\
			&\geq      \int _{t_1}^{t_2} g(s) {\rm d}B_s .
			\end{align*}
			Hence we conclude that $ \int _{t_1}^{t_2} g(s){\rm d}B_s  \leq -q_n + p_n$. By the definition of $A_n$, on event $A_n$ we can find a sequence of times $(t_{2k}, t_{2k+1})$ such that $t_{2k}< t_{2k+1} $ and $ \int _{t_{2k}}^{t_{2k+1}} g(s) {\rm d}B_s  \leq -q_n + p_n$. Now, if we define $G_{u,t} = \int_u ^t g(s) {\rm d}B_s  $, we see that $G_{1,t}$ converges a.s. since it is a martingale of bounded quadratic variation. Hence $\mathbb{P} (A_n) =0$, i.e. $\mathbb{P} (\liminf _{t \rightarrow \infty} X_t < 0    )=0$. 
			$\hfill \blacksquare$
			
				The next comparison result is intuitively obvious, however, it will be useful for comparing processes with different drifts. 
				\begin{prop}\label{prop:Sec4Dominance}
					Let $(C_t)_{t\geq 0}$ and $(D_t)_{t\geq 0}$ stochastic processes in the same Wiener space, that satisfy
					$\mathrm{d}C_t = f_1 ( C_t   ) {\rm d} t +     g(t) { \rm d}B_t    $, $\mathrm{d}D_t = f_2 ( D_t   ) {\rm d} t +     g(t) { \rm d}B_t    $ respectively, where $g,f_1,f_2$ are deterministic real valued functions. Assume that $f_1(x) > f_2(  x)$ for all$x \in \mathbb{R}$ and $C_{s_0}> D_{s_0}    $, then $C_{t} >  D_{t} \forall t \geq s_0   $ a.s..  
				\end{prop}
				\textbf{Proof:} Define $\tau =\inf\{ \tau> s_0 | C_\tau= D_\tau     \}$, and set $D_\tau=C_\tau=c$, for $\tau<\infty$. Now, from continuity of $f_1$, and $f_2$ we can find $\delta$ such that $f_1(x) >f_2(x), \forall x \in (  c-\delta   ,  c  ]  $. However,  for all $s$ we have
				$C_\tau - D_\tau -( C_s - D_s   )=-( C_s - D_s   )    =   \int_  s ^ \tau f_1(C_u) - f_2 ( D_u) {\rm d} u        $.  Thus, for $ s $ such that $  C_y, D_y \in (c-\delta,c    ) \forall y \in ( s,\tau   ) $ we have
				\begin{align*}
				0&>    -( C_s - D_s   ) \\
				&= \int_  s ^ \tau f_1(C_u) - f_2 ( D_u) {\rm d} u      \\
				&>0.
				\end{align*} 
				Therefore, $ \{\tau<\infty\}$ has zero probability.
				$\hfill \blacksquare$	
		
			In what follows, we will prove two important lemmas, corresponding to lemmas Lemma \ref{lem:ContLSup} and Lemma \ref{lem:ContLInf}, for the discrete case. We will assume that $X_n$ satisfies 
		\begin{equation}\label{eq:LemmasSup}
		X_{n+1} - X_n\geq\frac{f(X_n)}{n^{\gamma}} +  \frac{Y_{n+1}}{n^{\gamma}},\,\gamma \in (1/2,1) 
		\end{equation}
		where $f$ satisfies $\forall \epsilon>0,\exists c>0, f(x)\geq c,~x \in (-\infty,-\epsilon)$, and the $Y_n$ are defined similarly, as in \eqref{eq:genDiscNconv}.
		\begin{lemma}
			$\limsup X_n \geq 0$ a.s..
		\end{lemma}
		\textit{The proof is nearly identical as in the continuous case (Lemma \ref{lem:ContLSup})}
	
		\textbf{Proof:} We argue as before by contradiction: assume that $\limsup _{n \rightarrow \infty} X_n < 0   $ with positive probability. Pick $\delta >0   $ such that 
		$\limsup _{t \rightarrow \infty} X_n< -\delta$ with positive probability. Then, we can find a time $m$ such that $X_n \leq  -\delta$ $\forall n\geq m$. However, this has as an immediate consequence that 
		$\sum_{i=1}^{n} \frac{f(X_{i} )}{i^\gamma}  \rightarrow \infty   $. Now, the process $ G_{1,\infty} $ is a.s. finite, indeed since $ G_{1,n}$ is martingale it suffices to prove that $\sup _{n}E(G_{1,n}^2)<\infty$. Indeed
		\begin{align*}
		E(G_{1,n}^2) &=   E(\left( \sum_{i=m}^{n-1}\frac{Y_{i+1}}{i^{\gamma} }  \right )^2     )       
		=E(     \sum_{i=1}^{n-1}\frac{Y_{i+1}^2}{i^{2\gamma}}       ) \\
		&= \sum_{i=1}^{n}E(    \frac{Y_{i+1}^2}{i^{2\gamma}}       )
		\leq M^2\sum_{i=1}^{n-1}i^{-2\gamma}  \\
		&    \leq M^2\sum_{i=1}^{   \infty}i^{-2\gamma}.   \\
		\end{align*}
		Finally, the last two observations imply that $X_t \rightarrow \infty$, a contradiction.
		$\hfill \blacksquare$

		\begin{lemma} \label{discInf}
			$\liminf _{t \rightarrow \infty} X_t \geq 0  $ a.s..
		\end{lemma}
		\textit{The proof is nearly identical as in the continuous case (Lemma \ref{lem:ContLInf})}
		
		\textbf{Proof:} We argue by contradiction. Assume that $\liminf _{n \rightarrow \infty} X_n < 0   $ on a set of positive probability. Take an enumeration of the pair of  positive rationals $(q_m,p_m)$ such that $q_m>p_m$. Now, define $A_{m}= \{     X_n \leq -q_m ~\text{i.o.}, X_n\geq  -p_m ~ \text{i.o.}   \}  $. Since $ \limsup _{n \rightarrow \infty} X_n \geq 0  $, we have $  \bigcup_{n\geq 0} A_{m}   =  \{ \liminf _{n \rightarrow \infty} X_n < 0           \}  $. Choose $ n_1<n_2$ such that $ X_{n_1} \geq -p_m    $ and 
		$ X_{n_2} \leq -q_m    $. Therefore $X_{n_2} -X_{n_1} \leq -q_m + p_m     $; however,
		\begin{align*}
		X_{n_2} - X_{n_1} &= \sum_{i=n_1}^{n_2-1} \frac{f(X_i)}{i^\gamma}   + G_{n_1,n_2}    \\
		&\geq     G_{n_1,n_2}  
		\end{align*}
		Hence we conclude that $ G_{n_1,n_2}    \leq -q_m+ p_m$. By the definition of $A_m$ there is a sequence of times $(n_{2k}, n_{2k+1})$ such that $n_{2k}< n_{2k+1} $ and $  G_{n_1,n_2} \leq -q_m + p_m$. Hence $\mathbb{P} (A_m) =0$. Therefore $\mathbb{P} (\liminf _{n \rightarrow \infty} X_n < 0    )=0$. 
		$\hfill \blacksquare$

			\end{section}

			\section{Continuous model, simplest case}\label{sec:k=1}
			\subsection{Introduction}
			Let $L_t$ be defined by \eqref{eq:GenIntro}, for $f(x)=k|x|$ and $\gamma=1$. To simplify, we make a time change and consider $X_t := L_{e^t}$, and subsequently we obtain,
			\begin{align*}
			X_{t+{\rm d}t} - X_t &= L_{e^t+e^t {\rm d}t} -L_{e^t}\\
			&= k|L_{e^t}| {\rm d}t + e^{-t} (B_{t + e^t {\rm d}t} - B_{e^t}) \\
			&= k|X_t| {\rm d}t + e^{-\frac{t}{2}} {\rm d}B_t.
			\end{align*}
			Which will be the model we will study.
			We begin with some definitions.	
			\begin{equation} \label{eq:absval}
			{\rm d}X_t  = k |X_t| {\rm d}t + e^{-\frac{t}{2}}{\rm d} B_t.
			\end{equation}
			We introduce another SDE closely related to the previous one, which will be useful. 
			\begin{equation}\label{eq:noabsval}
			{\rm d}K_t  = k K_t {\rm d}t + e^{-\frac{t}{2}}{\rm d} B_t.
			\end{equation}
			It is easy to see that both of these SDEs admit unique strong solutions, for a reference see theorem $(11.2)$ in chapter $6$ in~\cite{Rogers87diffusions}. Therefore, we can construct $X_t,K_t$ in the classical Wiener space $(\Omega, \mathscr{F}, \mathbb{P})$. The solution for SDE~\eqref{eq:noabsval}, is given by $ K_t   =  e^{ k t} (e^{ -t_0k   } K_{t_0}    + \int _{t_0}^t e^{-s (k +\frac{1}{2}   )   }  {\rm d} B_s        )$. Indeed, substituting in $\eqref{eq:noabsval}$, and using It\^{o}'s formula we get 
			\begin{align*}  {\rm d} K_t   &= a'(t) (k_0 + \int _{t_0}^{t} b(s) {\rm d} B_s       ) + a(t)b(t){\rm d} B_t \\
			&= \frac{a'(t)}{a(t)} K_t+ a(t)b(t) {\rm d}B_t .  
			\end{align*} 	  
			Where
			$ a(t)= e^{ k (t-t_0)} $, and $ b(t) = e^{-t(\frac{1}{2}+k) +kt_0     }$. Therefore, $\frac{a'(t)}{a(t)} = k$, $ a(t)b(t)= e^{-\frac{t}{2}} $, so we conclude. 
			\begin{prop}
				Let $(X_t)_{t \geq t_0},\,(K_t)_{t \geq t_0}$ in the Wiener probability space $(\Omega,\mathscr{F},\mathbb{P})$ be the solutions of \eqref{eq:absval},\eqref{eq:noabsval} respectively.
				We start them at time $t_0$, $ X_{t_0}  \geq K_{t_0} \geq 0$. Then $ X_t \geq K_t, \, \forall t\geq t_0   $.
			\end{prop}
			It is a direct application of Proposition \ref{prop:Sec4Dominance}.
			$			\hfill \blacksquare$

			\subsection{Analysis of $X_t$ when $k>1/2$.}\label{sbec:k=1Nconv}
			Throughout this section, we assume that the parameter $k$ in \eqref{eq:absval} satisfies $k>1/2$.
			We start by stating the main result of this section, namely,		
			\begin{thm}\label{thm:maintheoremk=1}
			Let $(X_t)_{t\geq 1}$ the solution of \eqref{eq:absval}, then $X_t \rightarrow \infty$ a.s..
			\end{thm}
			We will prove the theorem at the end of the section. Now, we will show that $(X_t)_{t\geq 1}$ cannot stay negative for all times. This will be accomplished by a direct computation, after solving the SDE.
			\begin{prop}\label{prop:neverreturn}
				Assume that at time $s$, $X_{s}   < 0 $, then $X_t$ will reach $0$ with probability $1$, i.e. $\mathbb{P}( \sup_{u\geq s} X_u > 0        )=1$ 
			\end{prop}
			\textbf{Proof:} First, note that the solution of the SDE \eqref{eq:absval}, run from time $s$ with initial condition $X_{s}<0$ coincides with the solution of the SDE 
			$ {\rm d}X_t  = -k X_t {\rm d}t + e^{-\frac{t}{2}}{\rm d} B_t    $ before $X_t$ hits $0$. Formally, we define $\tau _0 = \inf \{ t| t \geq s,\,X_t=0    \}$. Using the same method when solving SDE \eqref{eq:noabsval}, we obtain $X_t =e^{-kt} (e^{ks} X_{s} + \int _{s}^{t} e^{u(k-\frac{1}{2}) } {\rm d}B_u      ) $, on $\{t <\tau _0 \}$. Set 
			$G_t= \int _{s}^{t} e^{u(k-\frac{1}{2}) } {\rm d}B_u      $, and calculate the quadratic variation of $G_t$, namely $\left \langle G_t \right \rangle =   ( e^{2t(k-\frac{1}{2}) }   -e^{2s(k-\frac{1}{2}) }       )/(2k-1)      $.     Next, we compute the probability of never returning to zero. 
			\begin{align*}
			\mathbb{P}( \tau =\infty      )&=	\mathbb{P}\left ( \sup _{  s <u <\infty } X_u \leq 0      \right )\\
			& =  	\mathbb{P}\left( \sup _{  s <u <\infty } G_u \leq-e^{ks} X_{s}        \right )\\ 
			&= 1- \mathbb{P}\left( \sup _{  s <u <\infty } G_u >-e^{ks} X_{s}   \right ) \\
			&=1-  \lim _{t\rightarrow \infty }\mathbb{P}\left( \sup _{  s <u <t } G_u >-e^{ks} X_{s}       \right  ) \\
			&= 1-  \lim _{t\rightarrow \infty }2 \mathbb{P}\left(  G_t >-e^{ks}  X_{s}\right ), \text{~ from\,the\,reflection\,principle} \\
			&= 1-    \lim _{t\rightarrow \infty }2 \mathbb{P}\left(  N\left (  0, \frac{   e^{2t(k-\frac{1}{2}) }   -e^{2s(k-\frac{1}{2}) } }{2k -1      } \right ) >-e^{ks} X_{s} \right ) \\
			&=0, \text{~since\,the\,variance\,is\,going\,to\,infinity}.
			\end{align*}
			$\hfill \blacksquare$
		
			We~will now prove an important lemma.
			\begin{lemma}\label{lemma2.2BecomesBfromRight} 
				There is a positive constant $c<1$, such that the set $ \{s| X_s \geq c e^{-s/2}  \}$ contains, with probability 1, an increasing sequence $t_n$ such that 
				$t_n\rightarrow \infty$ a.s., i.e. $\{  X_t \geq c e^{-s/2} \, \text{i.o.} \}$ happens a.s..	
			\end{lemma}	 
			Remark: To show the existence of the previously introduced times. We will see that the process $X_t$ starting at time $s$ such that $X_{s}\geq 0$ has enough variance to establish the existence of a future time $z$ such that 
				$X_z > ce^{ -z/2} $ with positive probability independent of time $s$.
			
			\textbf{Proof:} Since the process, when negative, eventually reaches zero with probability 1, a.s. there are times 
			$t_i$ such that $ t_{i+1} >  g(t_i)$, and $X_{t_i}\geq 0$, for any increasing function $g$. The function $g$ will be chosen shortly to ensure that the variance of $(X_t|\mathscr{F}_{t_i}  )$, is big enough. 
			
			Assume we start the SDE at time $t_i$ with initial condition $X_{t_i}\geq 0$. Then we see that $X_t \geq \int _{t_i}^t k|X_u| {\rm d }u + \int _{t_i}^t e^{-\frac{u}{2}} {\rm d}B_u   \geq \int _{t_i}^t e^{-\frac{u}{2}} {\rm d}B_u.$
			
			Set $G_t =\int _{t_i}^t e^{-\frac{u}{2}} {\rm d}B_u$. The quadratic variation of $G_t$, is $\langle G_t\rangle =  e^{-t_1} - e^{-t}$. Now, we observe that we can always choose $t$ big enough such that 
			$\langle G_t \rangle \geq c e^{-t_1}$.
		
			Then,
			\begin{align*}
			\mathbb{P} (  \sup _{t_i<u<t }X_t > e^{-t_1/2 } ) &\geq    \mathbb{P} (  \sup _{t_i<u<t }G_t > e^{-t_1/2 } )   \\
			&=     2\mathbb{P} ( G_t > e^{-t_1/2 } ) \\
			&\geq 2\mathbb{P} ( N(0,c e^{-t_1}) > e^{-t_1/2 } ) \\
			&= 2\mathbb{P} ( N(0, c) >1 ) > \gamma >0.
			\end{align*}
			Let $ g(x) =  \inf  \{  y| e^{-x} - e^{-y} \geq   c e^{-x}   \}$. Now, we can formally define the sequence of the stopping times. The first stopping time is $t_1 = \inf \{t |X_{t} \geq 0     \} $, then we define recursively 
			$t_{i+1} = \inf  \{ t | t> t_i, t> g(t_i), \, X_t \geq 0 \}$.
			
			To finish the proof, let us recall a suitable version of the Borel-Cantelli lemma (for a reference see Theorem $5.3.2$ in \cite{dur13}).
			\begin{lemma}{\label{Borel-Cantelli}}
				Let $\mathscr{F}_n$, $n \geq 0$ be a filtration with $\mathscr{F}_0 = \{ 0, \Omega \}$, and $A_n$, $n\geq 1$ a sequence of events with $A_n\in \mathscr{F}_n$. Then
				$$\{A_n \mathrm{\,i.o.} \} = \left \{ \sum_{n\geq 1} \mathbb{P}(A_n | \mathscr{F}_{n-1}   )  = \infty   \right  \}.$$
			\end{lemma} 
			We define the filtration $ \mathscr{F}_n = \mathscr{F}_{ t_n}$, for $n\geq 1$ and $\mathscr{F}_0 =\{ 0 , \Omega   \}$. Now, let $  A_n =\{ \exists t,\, t_{n-1} <t< t_{n} \, ,  \text{\,s.t.} X_t \geq  c e^{- t/2}       \}$. So, by definition $A_n \in \mathscr{F}_n$. 
			
			We find a lower bound for $ \mathbb{P} (A_n | \mathscr{F}_{n-1}    )$.
			\begin{align*}
			\mathbb{P} (A_n | \mathscr{F}_{n-1}    ) &\geq     \mathbb{P} ( \sup _{ t_{n-1} <u <  t_n } X_u > ce^{-t_{n-1} /2    } | \mathscr{F}_{n-1}    ) \\
			&\geq \mathbb{P} ( \sup _{ t_{n-1} <u < g(t_{ n-1})} X_u > ce^{-t_{n-1} /2    } | \mathscr{F}_{n-1}    ) \\
			&> \gamma. 	
			\end{align*}
			Hence $\sum_{n\geq 1} \mathbb{P}(A_n | \mathscr{F}_{n-1}   )  = \infty $ a.s. and we conclude. $\hfill \blacksquare$
			
			We know that for any initial value $X_{s}<0$ for any $s$, the process will eventually reach $0$. Before proving that $X_t\rightarrow \infty$, we will establish that the process will eventually stay positive. 
		\begin{lemma}\label{Lemma:2.2StaysPositive}
					Almost surely there is a time $u$ and a constant $\beta>0$ such that $X_t \geq \beta $ for all $t\geq u$,\, i.e. $\{\liminf_{t \rightarrow \infty} X_t >0   \}$ holds a.s..
				\end{lemma} 
				\textbf{Proof:} Indeed, if we start the process at time $s$ with initial condition $X_{s} \geq ce^{  \frac{-s}{2}}$, then the solution of \eqref{eq:absval}, before hitting $0$, is given by 
				$$X_t= e^{kt} \left (e^{-ks} X_{s} + \int _{s}^{t} e^{-u(k+\frac{1}{2}) } {\rm d}B_u     \right )\geq e^{kt} \left (ce^{-s(k +\frac{1}{2} ) } + \int _{s}^{t} e^{-u(k+\frac{1}{2}) } {\rm d}B_u    \right )  .$$
				Denote $G_t=\int _{s}^{t} e^{-s(k+\frac{1}{2}) } {\rm d}B_s  $. 
				We calculate its quadratic variation $$\langle G_t \rangle = \dfrac{  e^{-2tk-t   }  }{ -2k-1   }   +\dfrac{  e^{-2sk-s   }  }{ 2k+1   } .$$ 
				Taking $ t \rightarrow \infty $, shows 
				$\langle G_{\infty }\rangle=  \dfrac{  e^{-2sk-s}   }{2k+1     }$.
				Therefore, 
				\begin{align}
				\mathbb{P} ( \inf _{s \leq u <\infty } X_u > 
				\frac{c}{2}e^{   -\frac{s}{2}  }  )&=\mathbb{P} ( \inf _{s \leq u <\infty } e^{ku}(  ce^{-s(k +\frac{1}{2})} +   G_u ) >  \frac{c}{2}e^{   -\frac{s}{2}  }  )\nonumber\\
				&\geq \mathbb{P} ( \inf _{s \leq u <\infty } e^{ks}(  ce^{-s(k +\frac{1}{2})} +  G_u ) >  \frac{c}{2}e^{   -\frac{s}{2}  }  )\nonumber \\
				&=\mathbb{P} ( \inf _{s \leq u <\infty }   ce^{-s(k +\frac{1}{2})} +   G_u  >  \frac{c}{2}e^{  -s(k +\frac{1}{2})  }  )\nonumber	\\
				&=\mathbb{P} ( \inf _{s \leq u <\infty }   G_u  >  -\frac{c}{2}e^{  -s(k +\frac{1}{2}  )} )\label{eq:2.2InfBigger}	\\
				&= 1- \mathbb{P} ( \sup _{s \leq u <\infty }   G_u  >\frac{c}{2}e^{ - s(k +\frac{1}{2} ) }   )\nonumber  \\
				&= 1-2\lim _{t\rightarrow \infty}  \mathbb{P} ( G_t>-\frac{c}{2}e^{  -s(k +\frac{1}{2})  }   ),\text{~by\,the\,reflection\,principle} \nonumber\\
				&=1-2\mathbb{P} ( N(0, \dfrac{  e^{-s(2k+1})   }{2k+1     }) > \frac{c}{2}e^{  -s(k +\frac{1}{2} ) }   )\nonumber \\
				&=1-2\mathbb{P} ( N(0,\frac{1}{k+1}) > \frac{c}{2}  )>\delta>0. \nonumber
				\end{align}
				Now, as in the proof of Lemma \,\ref{lemma2.2BecomesBfromRight}, an application of Borel-Cantelli (Lemma \ref{Borel-Cantelli}), shows that $\{X_t\geq c e^{-\frac{t}{2}} ~\text{i.o.}   \}$ holds a.s.. Therefore, if we define $\tau_0=0$, and $\tau_{n+1} = \{t>\tau_n+1| X_t \geq c e^{-\frac{t}{2}}    \}$ we see that $\tau_n<\infty$ a.s., and $\tau_n \rightarrow \infty$ a.s.. Also, we define the corresponding filtration, namely $\mathscr{F}_n = \sigma( \tau_n)$.
			
				To show that $A=\{\liminf_{\rightarrow \infty} X_t\leq 0\}$ has probability zero, it suffices to argue that there is a $\delta$ such that $ \mathbb{P}( A|\mathscr{F}_n   )<1-\delta, $ a.s. for all $n\geq 1$. This is immediate from the previous calculation. Indeed,
				\begin{align*}
				\mathbb{P}( A|\mathscr{F}_n   ) & \leq  \mathbb{P} ( \inf _{ \tau_{n}\leq u <\infty } X_u > 
				\frac{c}{2}e^{   -\frac{\tau_n}{2}  } |\mathscr{F}_n )\\
				&<1-\delta.
				\end{align*}
				$\hfill \blacksquare$
			
				We can now show that $X_t \rightarrow \infty $ a.s.. 
				
				\textbf{Proof of Theorem \ref{thm:maintheoremk=1}:} Since, from Lemma \,\ref{Lemma:2.2StaysPositive}, $\liminf _{t \rightarrow \infty } X_t >0$ a.s., we deduce that 
				$ \int _{0}^\infty |X_u| { \rm d} u \rightarrow \infty$ a.s. 
				At the same time $ \limsup _{t \rightarrow \infty} \int _{0}^{t}  e^{-\frac{u}{2}} {\rm d } B_u < \infty$ a.s., hence $X_t\rightarrow \infty$ a.s.				
				$\hfill \blacksquare$

				\subsection{Analysis of $X_t$ when $k<1/2$.}\label{sbec:k=1C}
				As before, $(X_t)_{t\geq 1}$ is the solution of the stochastic differential equation
				$dX_t= k |X_t|\rm{d}t   +e^{-\frac{t}{2}}   \text{d}B_t $.
				 
				The behavior of $X_t$, when $k<1/2$ is different. The process in this regime can converge to $0$ with positive probability. More specifically, we have the following theorem:
				\begin{thm}\label{thm:contconvmain} 
					Let $(X_t)_{t \geq 1}$ solve \eqref{eq:absval} with $k<\frac{1}{2}$, and define $A= \{X_t\rightarrow 0 \}$, $B=\{X_t\rightarrow \infty     \}$. Then the following hold: 
					\begin{enumerate}[i)]
						\item  Let $A,B$ as before. Then $\mathbb{P}( A\cup B)=1$.
						
						\item Both $A$ and $B$ are non trivial i.e. $\mathbb{P}(A) >0$ and $\mathbb{P}(B) >0.$
						
						\item  On $\{X_t\geq 0~\mathrm{i.o.} \} $ we get $X_t\rightarrow \infty.$
					\end{enumerate}	
				\end{thm}
				We will prove the theorem at the end of this section. The next proposition shows that the process, starting from a negative value, will never cross $0$ with positive probability.				
				\begin{prop}\label{prop:statyingnegative}
					Assume that at time $s$, $X_{s}   < 0 $. Then $(X_t)_{t \geq 1}$ will hit $0$ with probability $\alpha$ where $0<\alpha<1$. 
				\end{prop}
				\textbf{Proof:} Define the stopping time $\tau_1  =\inf\{ t\geq s| X_t =0    \}$. As in Proposition \ref{prop:neverreturn}, the solution for $X_t$ started at time $s$ up till time $\tau_1$,is given by $X_t =e^{-kt} (e^{ks} X_{s} + \int _{s}^{t} e^{u(k-\frac{1}{2}) } {\rm d}B_u      ) $.
				\begin{align*}
				\mathbb{P}( \tau =\infty      )&=	\mathbb{P}( \sup _{  s <u <\infty } X_u \leq 0      )\\ 
				&= 1-\lim _{t\rightarrow \infty }2 \mathbb{P}(  N(  0, \frac{   e^{2t(k-\frac{1}{2}) }   -e^{2t(k-\frac{1}{2}) } }{2k -1      } ) >-e^{ks} X_{s}),\text{~as\,in\,Proposition\,}\ref{prop:neverreturn} \\
				&= 1-    2 \mathbb{P}(  N(  0,     -e^{2s(k-\frac{1}{2}) }       /(2k-1)      ) >-e^{ks} X_{s})=1-\alpha.
				\end{align*}
				Therefore $0<\alpha<1$.
				$\hfill \blacksquare$
				\begin{prop}\label{prop:comparisonk=1}
					Suppose
				 $(X_t)_{t\geq 1}$, $(Y_t)_{t \geq 1}$ solve \eqref{eq:absval}, with constants $k$ and $k_1$ respectively. Suppose, that $0<k_1<k<1/2$. If $X_s,Y_s \in(-\delta,0)$, then for all $\epsilon >0$ small enough, there is an event $A$ with positive probability, such that both $X_t,Y_t\in (-\delta-\epsilon,0),~\forall t>s$.
				\end{prop}
				\textbf{Proof:} Solving the SDE before it hits zero we find,
				$X_t =e^{-kt} (e^{ks} X_{s} + \int _{s}^{t} e^{u(k-\frac{1}{2}) } {\rm d}B_u      ) $ and $Y_t =e^{-k_1t} (e^{k_1s} Y_{s} + \int _{s}^{t} e^{u(k_1-\frac{1}{2}) } {\rm d}B_u      ) $.
				Let $\epsilon >0$, such that $\epsilon<\delta$. Since, the process $G_{t}=\int _{s}^{t} e^{u(k-\frac{1}{2}) } {\rm d}B_u    $ has finite quadratic variation, the event $A=\{G_{t}\in(-\epsilon,\epsilon) ~\forall t>s        \}$ has positive probability. Set $\tilde{G}_{t}=\int _{s}^{t} e^{u(k_1-\frac{1}{2}) } {\rm d}B_u   $, and define $N_t= G_te^{t(k_1-k)}$.
				Using It\^{o}'s formula, we find $\text{d}N_t= 
				e^{t(k-\frac{1}{2})} e^{(k_1-k) t}  \text{d}B_t+ (k_1-k)e^{(k_1-k) t}G_t\text{d}t$.
				Therefore, 
				$$G_te^{(k_1-k) t}= \tilde{G}_t+ \int_{s}^t e^{(k_1-k) u} G_u\text{d}u  .$$
				So, $$G_te^{(k_1-k) t}-\int_{s}^t (k_1-k)e^{(k_1-k) u} G_u\text{d}u= \tilde{G}_t  .$$
				Thus on $A$, we obtain the following inequalities
				$$-\epsilon e^{(k_1-k) t}+\epsilon( e^{(k_1-k) t}- e^{(k_1-k) s})\leq \tilde{G}_t \leq \epsilon e^{(k_1-k) t}-\epsilon( e^{(k_1-k) t}- e^{(k_1-k) s}) .  $$
				Simplifying, we obtain $ | \tilde{G}_t    | \leq \epsilon e^{(k_1-k) s}\leq \epsilon$.
				Now, we will estimate $X_t$ on $A$. We have the following upper bound,
				\begin{align*}
				X_t &= e^{-kt} (e^{ks} X_{s} + \int _{s}^{t} e^{u(k-\frac{1}{2}) } {\rm d}B_u      ) \\
				&\leq e^{-kt} (e^{ks} X_{s} +\epsilon )\\
				&<0 .
				\end{align*}
				and  lower bound
		\begin{align*}
		X_t &= e^{-kt} (e^{ks} X_{s} + \int _{s}^{t} e^{u(k-\frac{1}{2}) } {\rm d}B_u      ) \\
		&\geq e^{-kt} (-e^{ks} \delta +\epsilon )\\
		&\geq -\delta+\epsilon
		\end{align*}
		Doing similarly for $Y_t$, we conclude.
		$\hfill \blacksquare$
				\begin{prop}\label{prop:liminfk=1}
					Let $\epsilon >0$, then the event $   \{X_t <-\epsilon~ {\rm i.o.}    \}    $ has probability zero.
				\end{prop}
				\textbf{Proof:} This is a direct application of Lemma \ref{lem:ContLInf}.
				$\hfill \blacksquare$
				
				 \textbf{Proof of Theorem \ref{thm:contconvmain}:} \begin{enumerate}[i)] \item Define the events $N=\{ \exists s,\, \text{s.t.} X_t <0 \forall t\geq s       \}$, and $P=\{X_t \geq 0\, \text{i.o.}    \}$. Of course $N$ and 
					$P$ are disjoint and $\mathbb{P}(P \cup N   )=1$. To prove $\textit{i}$, we will show that 
					$N\subset \{ X_t \rightarrow 0   \}$ and $P= \{ X_t \rightarrow \infty    \}    $.
					
					Notice that Proposition \ref{prop:liminfk=1} implies that $ N \cap \left(            \bigcap _{k\geq 1}\{X_t <-1/k~ {\rm i.o.}    \}^{ \text{c}  }  \right) =N$ which exactly shows that $N \subset \{ X_t \rightarrow 0   \}$.
					
					To show that  $P= \{ X_t \rightarrow \infty    \}    $, we just need to make two observations. First, note that Lemma \ref{lemma2.2BecomesBfromRight}, actually proves something stronger, it shows that on $\{ X_t \geq 0\, \text{i.o.}   \}  $, $X_t \geq c e^{-\frac{t}{2}}$ holds infinitely often. Consequently, Lemma \ref{lemma2.2BecomesBfromRight} shows that on $ \{ X_t \geq c e^{-\frac{t}{2}}\, \text{i.o.}   \}$, $X_t\rightarrow \infty $. Therefore, $P= \{ X_t \rightarrow \infty    \}    $. Which concludes part $\text{i)}$.
					
					\item The fact that $\mathbb{P}(A)>0$, follows immediately from Proposition \ref{prop:statyingnegative}. Now, we will prove that $\mathbb{P}(B)>0$.
					Define the stopping time $\tau_0=\inf\{t |X_t=0     \}$. Also, define $Y(t,\omega)=1$ if $X_s\geq 0$ for all $s\geq t+1$. Observe that 
					$\{ Y_{\tau_0}=1,  \tau_0 < \infty   \}  \subset P $.
					Hence, using the strong Markov property
					\begin{align*}
					\mathbb{P}(Y_{\tau}=1,  \tau < \infty  ) &= \int_{0}^{\infty} \mathbb{ P} ( \tau=u    ) \mathbb{P}_0 ( X_t\geq 0, \forall t\geq 1 ){\rm d} u\\
					&\geq \int_{0}^{\infty} \mathbb{ P} ( \tau=u    ) \mathbb{P}_0 ( K_t\geq 0, \forall t\geq 1 ){\rm d} u\,\text{\,since\,}\,X_t\geq K_t\\
					&=\alpha \mathbb{P}_0 ( K_t\geq 0, \forall t\geq 1 )>0.
					\end{align*}

					\item Follows immediately from the proof of $\text{i)}$.
				\end{enumerate}
				$\hfill \blacksquare$

		\section{Analysis of ${ \rm d  }L_t = \frac{|L_t|^k}{t^{\gamma}} { \rm d}t + \frac{1}{t^\gamma}   {\rm d} B_t  $.}\label{sec:k>1}

		\subsection{Introduction}
		As in the previous section, to simplify matters, we will work with reparametrizing $L_t$. Set $\theta (t) = t^{ \frac{1}{1-\gamma}  }$, and let $ X_t= L_{ \theta (t)}$. To obtain the SDE that $X_t$ obeys, notice that $ {\rm d} B_{ \theta (t)  }  = \sqrt{\theta ' (t)} {\rm d} B_t $. Therefore 
		\begin{align*}
		{ \rm d  }X_t &= \frac{|X_t|^k}{\theta(t)^{\gamma}} \theta'(t){ \rm d}t + \frac{1}{\theta(t)^\gamma}   \sqrt{\theta ' (t)}{\rm d} B_t \\
		&= c_1|X_t|^k { \rm d}t + c_2t^{-\frac{\gamma}{1-\gamma}}   \sqrt{\theta ' (t)}{\rm d} B_t \\
		&=c_1|X_t|^k { \rm d}t + c_2t^{-\frac{\gamma}{2(1-\gamma)}}   {\rm d} B_t
		\end{align*}
		where $c_2^2=c_1= 1/(1-\gamma)$. By abusing the notation we set $ X_t = X_t/c_2$, which satisfies an SDE of the form 
		\begin{equation}\label{eq:mainwithScale}
		{ \rm d  }X_t =c|X_t|^k { \rm d}t + t^{-\frac{\gamma}{2(1-\gamma)}}  {\rm d} B_t ~k>1~\text{and}~\gamma\in(1/2,1) .
		\end{equation}
		Where $c$ is a positive constant. By a time scaling, we may assume that $X_t$ solves
		\begin{equation}\label{eq:main}
		{ \rm d  }X_t =|X_t|^k { \rm d}t + t^{-\frac{\gamma}{2(1-\gamma)}}  {\rm d} B_t,~k>1~\text{and}~\gamma\in(1/2,1) .
		\end{equation} 
		Notice, that the noise is scaled differently. However, it will be evident that only the order of the noise is relevant. 
		
		The solution of the SDE \eqref{eq:main}, when $X_t$ is positive, explodes in finite time. However, since we are interested in the behavior of $X_t$ when $X_t<M$ for a positive constant $M$, we may change the drift when $X_t$ surpasses the value $M$, which in turn it would imply that SDE \eqref{eq:main} admits strong solutions. One way to do this is by studying the SDE whose drift term is equal to $|x|^k$ when $x<M$ and $M$ when $x>M$. This SDE can be seen to admit strong solutions for infinite time. The reason is that this process $X_t$ is a.s. bounded from below, as the drift is positive. Also, $X_t$ cannot explode  to plus infinity in finite time since the drift is bounded from above when $X_t$ is positive. However, for simplicity, we will use the form as shown in \eqref{eq:main}.

		\subsection{Analysis of $X_t$ when $1/2 + 1/2k \geq \gamma$, $k>1$ and $\gamma \in (1/2,1)$}\label{sbec:ContNConv}
	    The main result of this section is the following theorem.
		\begin{thm}\label{thm:contnonconvnongen}
			Let $(X_t)_{t\geq 1}$, that solves $\eqref{eq:main}$. When  $1/2 + 1/2k >\gamma$,\,$X_t\rightarrow \infty$ a.s.
		\end{thm}
		We will see its proof at the end of the section. Now, we will prove an important proposition.		
		\begin{prop}
			The process $X_t$ gets close to the origin (from the left) infinitely often. More specifically, for some $\beta<0$, the event $\{X_t \geq \beta t^{\frac{1-2\gamma}{2(1-\gamma)}   } \text{i.o.}  \} $ has probability $1$.
		\end{prop}
		\textbf{Proof:}
		Define $h(t)= -t^{ \frac{1}{1-k}}$, and notice that $h'(t) = -\dfrac{1}{1-k}  |h(t)|  ^k$. Then define $Z_t = -\dfrac{X_t}{h(t) }     $. Recall that since $h(t)$ is a continuous function, the covariance $ \langle h(t) , Z_t \rangle =0$. Now, we will find the SDE that $Z_t$ satisfies when $X_t$ solves \,\eqref{eq:mainwithScale} i.e. the drift of $X_t$ has a scaling factor. And, right after the computation ends, we will set $c=1$ again. 
		We use It\^{o}'s formula, and obtain 
		$$ {\rm d} Z_t =  -\frac{1}{h(t)} {\rm d} X_t  +   X_t {   \rm d } \left ( -\frac{1}{h(t) }  \right ) .$$ 
		Thus,
		\begin{align*}
		Z_t - Z_s &= \int _{ s  }^{t}   -\frac{1}{h(u)} c|X_u| ^ k { \rm d } u + \int _{ s  }^{t}   -\frac{1}{h(u)}u^{-\frac{\gamma}{2(1-\gamma)}}   { \rm d }B_u + \int _{ s  }^{t}  X_u\frac{h'(u)}{ h(u) ^2} {\rm d }u       \\  
		&=  \int _{ s  }^{t}  X_u\frac{h'(u)}{ h(u) ^2}   -\frac{1}{h(u)} c|X_u| ^ k { \rm d } u + \int _{ s  }^{t}   -\frac{1}{h(u)}u^{-\frac{\gamma}{2(1-\gamma)}}   { \rm d }B_u   \\
		&=\int _{ s  }^{t}c\dfrac{X_u}{h(u) }  \left (\frac{h'(u)}{c h(u) }   - \frac{|X_u| ^k }{X_u}\right ) { \rm d } u + \int _{ s  }^{t}   -\frac{1}{h(u)}u^{-\frac{\gamma}{2(1-\gamma)}}   { \rm d }B_u   \\
		&=\int _{ s  }^{t}c\dfrac{X_u}{h(u) }  \left (  \frac{1}{c(k-1)}\frac{|h(t)|^k }{h(t)}    -\frac{|X_u| ^k }{X_u}\right ) { \rm d } u + \int _{ s  }^{t}   -\frac{1}{h(u)}u^{-\frac{\gamma}{2(1-\gamma)}}   { \rm d }B_u.  
		\end{align*}
		This gives an SDE for $Z_t$ in terms of $X_t$ and $h(t)$
		\begin{equation}\label{eq:SDEunified}
		Z_t - Z_s= \int _{ s  }^{t}c\dfrac{X_u}{h(u) }  \left ( C\frac{|h(t)|^k }{h(t)}   -\frac{|X_u| ^k }{X_u}\right ) { \rm d } u + \int _{ s  }^{t}   -\frac{1}{h(u)}u^{-\frac{\gamma}{2(1-\gamma)}}   { \rm d }B_u, 
		\end{equation}
		where $C(c) =  \frac{1}{c(k-1)}$ is a function of $c$, and whenever there is no ambiguity what the argument is we will just call it $C$. We set $c=1$, and continue with the proof of the Proposition.
		Set ${G'}_t =\int _{ s  }^{t}   -\frac{1}{h(u)}u^{-\frac{\gamma}{2(1-\gamma)}}   { \rm d }B_u $. We calculate its quadratic variation at time $t$, $\langle{G'}_t\rangle =\int _{ s  }^{t}   \frac{1}{h(u)^2}u^{-\frac{\gamma}{(1-\gamma)}}   { \rm d }u $
		Notice from the definition of $h(t)$ we have $h(t) = \Theta (t^{\frac{1}{1-k}})$ so
		$h(t)^{-1} = \Theta (t^{\frac{1}{k-1}})$.
		So the integrand  $-\frac{1}{h(u)^2}u^{-\frac{\gamma}{(1-\gamma)}} = \Theta (u^{\frac{2}{k-1}-\frac{\gamma}{1-\gamma   }     })  $. Therefore, $ \langle {G'}_{\infty}\rangle=\infty$ when $$\frac{2}{k-1}-\frac{\gamma}{1-\gamma   }    \geq-1 \iff \frac{2}{k-1}+\frac{1}{\gamma-1   }     \geq -2.$$
		Now, from the restrictions on $k$ we obtain
		\begin{align*}
		\frac{1}{2} +\frac{1}{2k} \geq \gamma &\iff \frac{1-k}{2k} \leq\gamma -1 \iff \frac{k-1}{k} \geq 2(1-\gamma )\\&\iff
		k-1 \geq \frac{2(1-\gamma )   }{2\gamma -1} \iff
		\frac{1}{k-1}  \leq \frac{ 2\gamma -1  }{2(1-\gamma )} 
		\end{align*} 
		Therefore, we get that 
		$\dfrac{    1}{1-\gamma} \geq \dfrac{  2k   }{k-1}$, so
		$  \dfrac{2}{k-1}+\dfrac{1}{\gamma-1   } \geq \dfrac{2}{k-1}+ \dfrac{2k}{1-k}=-2.$
	
	This is also gives
		$   - t^{   \frac{1}{1-k}     }\geq  - t^{ \frac{1-2\gamma     }{2(1-\gamma)}      }$.
	
		First, we will prove that $  \{  X_t \geq C' h(t),\,\text{i.o.}    \}$ a.s., where $C'> C^{\frac{1}{k-1}}$. To do so, we will argue by contradiction. Assume that $A=\{  \exists\,s, X_t<C'\cdot h(t) \forall t>s     \}$ has positive measure. Take $\omega \in A$, and find $s(\omega)$ such that $ X_t<C'\cdot h(t)$ for all $t>s$. Notice, that this implies that $Z_t <-C'$ for $t>s$.  Take $u>s$, since $ \frac{|x|^k}{x}$ is increasing we see that   $\frac{ |X_u|^k    }{X_u}< {C'}^{k-1} \frac{ |h(u)|^k    }{h(u )}     < C \frac{ |h(u)|^k    }{h(u )}  $. This in turn gives $  C \frac{ C|h(u)|^k    }{h(u )}  -\frac{ |X_u|^k    }{X_u}>0$. Therefore 
		$\int _{ s  }^{t}\dfrac{X_u}{h(u) }  \left (C\frac{ |h(u)|^k    }{h(u )}  -\frac{ |X_u|^k    }{X_u}\right ) { \rm d } u >0$ for all $t>0$. However, since $G'^w_t$ for any fixed $w$ has infinite quadratic variation, we may find $-G'^s_t >-Z_s$. Now, from~\eqref{eq:SDEunified} we get
		\begin{align*}
		Z_t &= \int _{ s  }^{t}\dfrac{X_u}{h(u) }  \left (C\frac{ |h(u)|^k    }{h(u )}  -\frac{ |X_u|^k    }{X_u}\right ) { \rm d } u + Z_t   -G'^s_t\\
		&>0.
		\end{align*}
		This contradicts the fact that $Z_t<-1$. Therefore $ \{X_t> C'h(t), \text{i.o.}     \}$ a.s.
		
		Finally, since  $ -t^{\frac{1}{k-1}} \geq  -t^{ \frac{1-2\gamma}{ 2(1-\gamma)}}$ we conclude that there exists a constant $\beta<0$ such that $\{ X_t \geq \beta t^{ \frac{1-2\gamma}{ 2(1-\gamma)}    }  \text{\,i.o.}   \}  $ holds a.s.
		$\hfill \blacksquare$
		\begin{cor}\label{Cor3.3}
			The event   $\{   X_t \geq 0\,\text{i.o.}   \}$ holds a.s.
		\end{cor}
		\textbf{Proof:} Set 
		$G_{s,t}=\int _s^t u^{-\frac{\gamma}{2(1-\gamma)}}  {\rm d} B_u .$ Notice the lower bound $ X_t-X_s \geq G_{t,s}$. Now, observe that 
		$  \langle G_{s,\infty}\rangle = \Theta (s^{ \frac{1-2\gamma}{ (1-\gamma)}    }      )$, and so a similar calculation as in the beginning of Lemma \,\ref{lemma2.2BecomesBfromRight} yields that for any $c>0$,
		$\mathbb{P}(\sup _{s<t<\infty}X_t-X_s \geq   c s^{ \frac{1-2\gamma}{ 2(1-\gamma)}    }  )> \delta$, for some constant $\delta>0$. Therefore, for $c>\beta$(as in Proposition 3.4) we see that given $\mathscr{F}_\tau$, where $\tau$ is a stopping time such that $X_\tau \geq -\beta \tau^{\frac{1}{1-k}}$, the probability that $X_t\geq 0$ for some $t>\tau$ is bigger than $\delta$. Now, using Lemma \ref{Borel-Cantelli} (Borel-Cantelli) as in Lemma \ref{lemma2.2BecomesBfromRight} we conclude.
		$\hfill \blacksquare$
		
		\textbf{Proof of Theorem \ref{thm:contnonconvnongen}:} First, recall that $  \langle G_{s,\infty}\rangle = \Theta (s^{ \frac{1-2\gamma}{ (1-\gamma)}    }      )$, and by a similar calculation as done in \,\eqref{eq:2.2InfBigger} in Lemma \,\ref{Lemma:2.2StaysPositive}, we obtain that $ \mathbb{P} ( \inf _{ s<t<\infty}G_{s,t} > -\dfrac{m}{2}  t^{ \frac{1-2\gamma}{ 2(1-\gamma)}    }      )>\delta>0  $ for every $m>0$.  Fix a time $s$, and notice that $ X_t-X_s \geq G_{s,t}$. Thus, observe that whenever $X_s \geq m  s^{ \frac{1-2\gamma}{ 2(1-\gamma)}    }  $, we see that $\mathbb{P}   ( \inf _ { s<   t< \infty   } X_t >  \dfrac{m}{2}  s^{ \frac{1-2\gamma}{ 2(1-\gamma)}    }      )   >\delta >0$. 
		
		From the proof of the Corollary \ref{Cor3.3}, we have that $\mathbb{P}(\sup _{s<t<\infty}X_t-X_s \geq   c s^{ \frac{1-2\gamma}{ 2(1-\gamma)}    }  )> \delta$, and since $\{X_t\geq 0 ~\text{i.o.} \} $ has probability $1$, Borel-Cantelli (Lemma \ref{Borel-Cantelli}) implies that the event $  \{ X_t \geq m  t^{ \frac{1-2\gamma}{ 2(1-\gamma)}    } \text{i.o.}    \}  $ holds a.s. Define $\tau_{n+1}=\inf \{t>\tau_n+1| X_t \geq m  t^{ \frac{1-2\gamma}{ 2(1-\gamma)}    }    \}$, and observe that $\tau_n<\infty$, a.s. From the last paragraph's observation $\mathbb{P}   ( \inf _ { \tau_n<   t< \infty   } X_t >  \dfrac{m}{2}  {\tau_n}^{ \frac{1-2\gamma}{ 2(1-\gamma)}    }  |\mathscr{F}_{\tau_n}    )   >\delta >0$, therefore similarly as in Lemma \ref{Lemma:2.2StaysPositive},  $\liminf_{t\rightarrow \infty}X_t>0     $ a.s. From here, we see that the drift term goes to infinity a.s., and the noise is finite a.s. Therefore, $X_t \rightarrow \infty$ a.s.

		$\hfill \blacksquare$

		\subsection{Analysis of $X_t$ when $\frac{1}{2} + \frac{1}{2k} < \gamma $, and $k>1$      }\label{sbec:contConv}
		We start by changing the form of \eqref{eq:SDEunified} by rewriting it in terms of $Z_t$. For the purposes of this section, we find the form of the SDE before $X_t$ has reached $0$, for an initial condition $X_s<0$; 
		\begin{equation}\label{eq:SDEmerged}
		Z_t - Z_s= \int _{ s  }^{t}c|h(u)|^{k-1} Z_u  \left ( C-   (-Z_u) ^ {k-1}\right ) { \rm d } u + \int _{ s  }^{t}   -\frac{1}{h(u)}u^{-\frac{\gamma}{2(1-\gamma)}}   { \rm d }B_u 
		\end{equation} 
		where recall $C=\frac{1}{c(k-1)}$.
		We now state the main theorem of this section, which we will prove at the end.
		\begin{thm}\label{thm:contconvnongen}
			The process $(X_t)_{t\geq 1}$ the solution of \eqref{eq:mainwithScale}, converges to zero with positive probability, when $X_1<0$. 
		\end{thm}
		Before proving the theorem we will need the following proposition.
		\begin{prop}\label{prop:Main3.3}
			Assume that at time $s$, $X_s<0$, and $Z_s>- \left (\frac{C}{k}   \right )^{\frac{1}{k-1}}$. Then the process with positive probability never returns to $0$.     
		\end{prop}
		\textbf{Proof:} 
		The condition $1/2 + 1/2k <\gamma$ as it has already been shown in previous section \ref{sbec:ContNConv}, implies that 
		$  \langle G_{\infty} \rangle <\infty$. 
		Now, we will prove a lemma needed for the proof. 
		\begin{lemma}
			Assume that $Z_s>-\left (\frac{C}{k}   \right )^{\frac{1}{k-1}}$. And let $A= \{G_t \in (-\epsilon, \epsilon    )\forall t \in (s, s+\delta), \text{and}\,\, G_t \in (-2\epsilon   ,-\frac{9}{10}\epsilon    ) \forall t \in  (s+\delta,\infty)\}$. Then:
			\begin{enumerate}\label{lemma:Main3.3}
				\item $\mathbb{P}(   A) >0,\,\forall \epsilon,\delta >0$.
				
				\item For all $\epsilon>0$ small enough,  there is $\delta>0$ such that $Z_t<-\dfrac{5\epsilon}{3},\,\, \forall t \in  (s, s+\delta)$ on $A$.
			\end{enumerate}
		\end{lemma}
		\textbf{Proof:} 
		\begin{enumerate}
		\item The first is immediate since $G_t$ has finite quadratic variation.
	
		\item The first restriction on $\epsilon$ so that $Z_s<-3\epsilon$. Next, we begin by defining $f_1$ and $f_2$ on $(s, s+\delta)$ satisfying  
		\begin{equation}\label{eq:IntIn3.3}
		f'(x) = c|h(x)|^{k-1} f(x) (C-(-f(x)^{k-1})      ) 
		\end{equation}
		with initial conditions $ -\left (\frac{C}{k}   \right )^{\frac{1}{k-1}} < Z_s +\epsilon< f_1(s ) <-\dfrac{5\epsilon}{3} $, and $ -\left (\frac{C}{k}   \right )^{\frac{1}{k-1}}<f_2(s )  < Z_s-\epsilon   $. 
		
		Also,  we define the function $ q(x) = x(C- (-x)^{k-1}    ) $, whose derivative is $ q'(x) = C-k(-x)^{k-1}$, which implies that $q(x)$ is increasing on $ (      \left (-\frac{C}{k}   \right )^{\frac{1}{k-1}} ,0     )   $. This function will be important later. We should also note, that $f$ is decreasing in intervals where $f(x) \in (   - \left (\frac{C}{k}   \right )^{\frac{1}{k-1}} ,0     )$, since there $f'(x)<0$.
		
		We can pick the $\delta>0$, such that $f_2(t)>- \left (\frac{C}{k}   \right )^{\frac{1}{k-1}}$    $,\,\forall t \in ( s, s+\delta ) $. We will show that $ Z_t>f_2(t)$ on $( s, s+\delta )   $ by contradiction.  Using \eqref{eq:SDEmerged}, we get that 
		\begin{equation}\label{eq:SDEestimate   }
		Z_t- Z_s =   \int _{ s  }^{t}c|h(u)|^{k-1} Z_u  \left (  C-   (-Z_u )^ {k-1}\right ) { \rm d } u + g(t),
		\end{equation}
		where $g(t)$ is a continuous function such that 
		$ \sup_{t\in ( s, s+\delta ) } |g(t)| \leq \epsilon. $
		Assume that $f_2, Z$ become equal at some point, and choose $t$ to be the first time. Using the integral form of \eqref{eq:IntIn3.3}, and subtracting it from \eqref{eq:SDEmerged}, we get 
		\begin{align*}
		0=Z_t- f_2(t) &= Z_s- f_2(s) +  \int _{ s  }^{t} c|h(u)|^{k-1} Z_u  \left ( C-  (-Z_u )^ {k-1}\right ) -c|h(u)|^{k-1} f_2(u)  \left ( C-   (-f_2(u)) ^ {k-1}\right )  { \rm d } u + g(t)  \\
		&=  Z_s+g(t) - f_2(s)   + (t-s)(c|h(\xi)|^{k-1} Z_{\xi}  \left ( C-  ( -Z_{\xi}) ^ {k-1}\right ) -c|h(\xi)|^{k-1} f_2(\xi)  \left ( C-  (-f_2(\xi)) ^ {k-1}\right ) ) \\
		&>(t-s)(c|h(\xi)|^{k-1} Z_{\xi}  \left ( C-   (-Z_{\xi}) ^ {k-1}\right ) -c|h(\xi)|^{k-1} f_2(\xi)  \left ( C-   (-f_2(\xi)) ^ {k-1}\right )),
		\end{align*}
		where in the last line we used that $Z_s+g(t) - f_2(s)> 0$.
		Since $\xi <t$, we have that $ Z_{\xi}>f_2(\xi) > -\left (\frac{C}{k}   \right )^{\frac{1}{k-1}}   $, and consequently  
		$ q(Z_{\xi} )    >   q(  f_2(\xi)   )   $, so
		$$  |h(\xi)|^{k-1} q(Z_{\xi} )    >   |h(\xi)|^{k-1} q(  f_2(\xi)   )    .$$ Therefore, 
		$$0 <  c|h(\xi)|^{k-1} Z_{\xi}  \left ( C-   (-Z_{\xi}) ^ {k-1}\right ) -c|h(\xi)|^{k-1} f_2(\xi)  \left ( C-   (-f_2(\xi)) ^ {k-1}\right ),    $$
		which gives a contradiction.
		
		To complete the proof, notice it suffices to show that $ f_1(t)> Z_t$ on $ (s, s+\delta)$, and to do so we will argue by contradiction again. As before, take $t$ to be the first time such that $f(\cdot)$ and $Z_t$ are equal.
	
		Using the integral form of \eqref{eq:IntIn3.3} and subtracting from \eqref{eq:SDEestimate   } we get 
		\begin{align*}
		0=Z_t- f_1(t) &= Z_s- f_1(s) +  \int _{ s  }^{t}c|h(u)|^{k-1} Z_u  \left ( C-  (- Z_u) ^ {k-1}\right ) -c|h(u)|^{k-1} f_1(u)  \left ( C-   (-f_1(u)) ^ {k-1}\right )  { \rm d } u + g(t)  \\
		&=  Z_s+g(t) - f_1(s)   +(t-s)(c|h(\xi)|^{k-1} Z_{\xi}  \left ( C-  (-Z_{\xi}) ^ {k-1}\right ) - c|h(\xi)|^{k-1} f_1(\xi)  \left ( C-  (- f_1(\xi)) ^ {k-1}\right ))  \\
		&<(t-s)(c|h(\xi)|^{k-1} Z_{\xi}  \left (  C-  (-Z_{\xi}) ^ {k-1}\right ) - c|h(\xi)|^{k-1} f_1(\xi)  \left ( C-   (-f_1(\xi)) ^ {k-1}\right )), 
		\end{align*}
		where in the last line we used that $Z_s+g(t) - f_1(s)< 0$.
		Since, $\xi <t$ we have that $-\left (\frac{C}{k}   \right )^{\frac{1}{k-1}}< Z_{\xi}<f_1(\xi)     $, consequently
		$ q(Z_{\xi} )   <  q(  f_1(\xi)   )   $, and by multiplying with $ |h(\xi)|^{k-1}   $ we obtain 
		$  |h(\xi)|^{k-1} q(Z_{\xi} )    <   |h(\xi)|^{k-1} q(  f_1(\xi)   )    $. Therefore, 
		$$ |h(\xi)^{k-1}| Z_{\xi}  \left ( 1-   (-Z_{\xi}) ^ {k-1}\right ) -|h(\xi)|^{k-1} f_1(\xi)  \left ( 1-   (-f_1(\xi)) ^ {k-1}\right )<0,    $$
		again a contradiction. 
		\end{enumerate}
		$\hfill \blacksquare$     
		
		We resume now to the proof of Proposition \ref{prop:Main3.3}.
		On $A$, using \eqref{eq:SDEunified}, we get the following upper and lower bounds for all $t\geq s+\delta$ 
		\begin{equation} \label{eq:UpperBound} 
		-\frac{X_t}{h(t)} \leq  -\frac{X_s}{h(s)} +\int _{ s  }^{t}c\dfrac{X_u}{h(u) }  \left ( C\frac{|h(t)|^k }{h(t)}   -\frac{|X_u| ^k }{X_u}\right ) { \rm d } u -\frac{9}{10}\epsilon     
		\end{equation}
		\begin{equation} \label{eq:LowerBound} 
		-\frac{X_t}{h(t)} \geq  -\frac{X_s}{h(s)} +\int _{ s  }^{t}c\dfrac{X_u}{h(u) } \left ( C\frac{|h(t)|^k }{h(t)}   -\frac{|X_u| ^k }{X_u}\right ) { \rm d } u -2\epsilon     
		\end{equation} 
		
		\textbf{Claim:} on $A$, $X_t < 0,$ for all $t>s$.
		
		\textbf{Proof:} We will argue by contradiction. Assume that $\mathbb{P}(\{  \tau_0 <\infty \}\cap A) >0$. We choose $\epsilon$, such that $ \frac{3\epsilon}{2}<C^{\frac{1}{k-1} }     $. Now, define $\tau_{l} = \sup \{  t\leq \tau_0 |-\dfrac{X_t}{h(t) } = -\frac{3\epsilon}{2}       \}    $ and notice that Lemma \ref{lemma:Main3.3}, implies that $ \tau _{  l \epsilon }>s+\delta$, since $Z_t< -\dfrac{5 \epsilon}{3}$, on
		$(s,s+\delta)$. Also, on $\{  \tau_0 <\infty \}\cap A$ we have $\tau _{ l}<\infty$.   Then from $\eqref{eq:LowerBound}$
		we see that 
		$$ \int _{ s  }^{ \tau _{  l}}c\dfrac{X_u}{h(u) } \left (C \frac{|h(t)|^k }{h(t)}   -\frac{|X_u| ^k }{X_u}\right ) { \rm d } u \leq \frac{X_s}{h(s)} +\frac{\epsilon}{2}  .$$
		Therefore, 
		\begin{equation} \label{eq:NozeroBound} 
		-\frac{X_s}{h(s)} +\int _{ s  }^{ \tau _{  l}}c\dfrac{X_u}{h(u) } \left ( C\frac{|h(t)|^k }{h(t)}   -\frac{|X_u| ^k }{X_u}\right ) { \rm d } u -\frac{9}{10}\epsilon       \leq -\frac{2\epsilon}{5}     .
		\end{equation}
		Now, notice that $ X_t > \frac{3}{2} \epsilon  h(t), \, \forall t \in ( \tau _{l},\tau_0  )$, so if $w \in ( \tau _{l},\tau_0  ),$ we get  $   C\frac{|h(w)|^k }{h(w)}   -\frac{|X_w| ^k }{X_w} < C\frac{|h(w)|^k }{h(w)}- C\frac{|h(w)|^k }{h(w)}  =0   $ and of course 
		$ \frac{X_w}{ h(w)}  >0$. So, we conclude that 
		\begin{equation} \label{eq:FinalBound}
		\int _{ \tau_l }^{\tau_0}c\dfrac{X_u}{h(u) }  \left (C \frac{|h(t)|^k }{h(t)}   -\frac{|X_u| ^k }{X_u}\right ) { \rm d } u <0.
		\end{equation}
		Combining $\eqref{eq:NozeroBound}$, and $\eqref{eq:FinalBound}$, we get that 
		\begin{align*}
		0=-\frac{X_{\tau_0}}{h(\tau_0)} &\leq  -\frac{X_s}{h(s)} +\int _{ s  }^{\tau_0}c\dfrac{X_u}{h(u) } \left (C \frac{|h(t)|^k }{h(t)}   -\frac{|X_u| ^k }{X_u}\right ) { \rm d } u -\frac{9}{10}\epsilon  \\
		&=-\frac{X_s}{h(s)} +\int _{ s  }^{\tau_l}c\dfrac{X_u}{h(u) }  \left ( C\frac{|h(t)|^k }{h(t)}   -\frac{|X_u| ^k }{X_u}\right ) { \rm d } u -\frac{9}{10}\epsilon+\int _{ \tau_l  }^{\tau_0}c\dfrac{X_u}{h(u) } \left ( C\frac{|h(t)|^k }{h(t)}   -\frac{|X_u| ^k }{X_u}\right ) { \rm d } u \\
		&\leq -\frac{2\epsilon}{5},
		\end{align*}
		a contradiction. 
		$\hfill \blacksquare$
		
		We have developed all the tools necessary, to prove the theorem.
		
		 \textbf{Proof of Theorem \ref{thm:contconvnongen}:} Define a stopping time $\sigma=\inf \{ t  | Z_t>-\left (\frac{C}{k}   \right )^{\frac{1}{k-1}}        \}$. If the event $\{\sigma<\infty\}$ has positive probability, then Proposition \ref{prop:Main3.3} implies that $X_t$ converges to zero with positive probability. Indeed, remember from Lemma \ref{lem:ContLInf} that $\liminf_{t\rightarrow \infty}X_t \geq 0$ a.s. Therefore, since on $A$ we have $\limsup_{t \rightarrow \infty }X_t \leq 0$, we deduce $\lim_{t\rightarrow \infty }X_t=0$. To finish the proof, it suffices to show that when $\{\sigma<\infty\}$ has zero probability then $X_t\rightarrow 0$ with positive probability. This is easy to see since $\mathbb{P}(\{\sigma<\infty\})=0$ implies that $Z_t$, never hits zero, therefore  $\limsup_{t \rightarrow \infty }X_t \leq 0$ on $\{\sigma<\infty \}$.
		$\hfill \blacksquare$
		
		We now prove a proposition that will be used in the next section. 
		\begin{prop}\label{prop:StayingClose}
			Take the event $A$, such that Lemma \ref{lemma:Main3.3} holds, where $\epsilon <\left (\frac{C}{k} \right)   ^{\frac{1}{k-1}}$. Then, $X_t$ on $A$ stays within a region of the origin. More specifically, $Z_t> -2\left (\frac{C}{k} \right)   ^{\frac{1}{k-1}}$.	
		\end{prop}
		\textbf{Proof:} Let $\tau_C = \inf \{t>s|   Z_t=    - 2 \left (\frac{C}{k} \right)   ^{\frac{1}{k-1}}        \}$, and define $\sigma= \sup \{\tau_C >t>s|Z_t=  -  \left (\frac{C}{k} \right)   ^{\frac{1}{k-1}}         \}$. We will show that $\tau_C=\infty$ a.s. We assume otherwise, and reach a contradiction. From the proof of Lemma \ref{lemma:Main3.3}, we know that $\tau_C>s+\delta$. Therefore,
		\begin{align*}
		Z_{ \tau_C}  &\geq  Z_s+  \int _{ s  }^{t}c|h(u)|^{k-1} Z_u  \left ( C-   (-Z_u) ^ {k-1}\right ) { \rm d } u    -2\epsilon              \\
		&\geq Z_{\sigma}  -Z_{\sigma}+Z_s+  \int _{ s  }^{t}c|h(u)|^{k-1} Z_u  \left ( C-   (-Z_u) ^ {k-1}\right ) { \rm d } u   -2\epsilon  \\
		&\geq Z_{\sigma}   +\frac{9\epsilon}{10}   -2\epsilon >-2\left (\frac{C}{k} \right)   ^{\frac{1}{k-1}},
		\end{align*}
		the desired contradiction.
		$\hfill \blacksquare$

		\begin{section}{ Analysis of ${ \rm d  }L_t = \frac{f(L_t)}{t^{\gamma}} { \rm d}t + \frac{1}{t^\gamma}   {\rm d} B_t  $    .           }\label{sec:ContGenCase}
				For this section, we assume that $f$ is globally Lipschitz.
			 For $f$ as before, we define 
			\begin{equation}\label{eq:genSDE}
			{ \rm d  }L_t = \frac{f(L_t)}{t^{\gamma}} { \rm d}t + \frac{1}{t^{\gamma}}  {\rm d} B_t, \gamma \in (\frac{1}{2},1]
			\end{equation}
			By our assumptions on $f$, the SDE \eqref{eq:genSDE} admits strong solutions.
			Also, we define a more general SDE, namely
			\begin{equation}\label{eq:genSDEgennoise}
			{ \rm d  }X_t = f(X_t) { \rm d}t + g(t)   {\rm d} B_t
			\end{equation}
			where $g:\mathbb{R}_{\geq 0}\rightarrow \mathbb{R}_{>0}$ is continuous, and $T=\int_{0}^\infty g^2(t)\text{d}t$ is possibly infinite.	
			\begin{prop}\label{prop:hittingeverypoint}
				Let $(X_t)_{t\geq 1}$ be a solution of \eqref{eq:genSDEgennoise}. Then for every $t,c>0$,  and $x\in \mathbb{R}$,  $\mathbb{P}( X_t \in(x-c ,x+c   )        )>0$.
			\end{prop}
			\textbf{Proof:} Firstly, we change time. Let $\xi(t)= \int _{0}^{t} g^2(t)\text{d}t $, and define $\tilde{X}_t =X_{\xi^{-1}(t)}$. Then, 
			\begin{equation}\label{eq:hittingeverytime}
			\text{d} \tilde{X}_t= \frac{f(\tilde{X}_t   )}{g^2(\xi^{-1}(t)   )    }\text{d}t       +\text{d}B_t     
			\end{equation}
			This gives a well defined SDE whose solution is defined on $[0,T']$ for  $T'\in \mathbb{R},$ $\xi(t)<T'<T$. The path space measure of $\tilde{X}_t$ is mutually absolutely continuous to the one induced from the Brownian motion. Since the Brownian motion satisfies the property described in the proposition, so does $X_t$.   
			$\hfill \blacksquare$ 
									
			We give the proof of Theorems \ref{thm:gennonconv} and \ref{thm:genconv}. For the proofs, we use that the theorems hold if and only if they hold for their corresponding  reparameterizations. 
			
			\textbf{Proof of Theorem \ref{thm:gennonconv} part \ref{enum:nonconv1} \& \ref{enum:nonconv2}:} Both parts can be proved simultaneously. Let $\tau = \inf \{ t|X_t \in (   -\epsilon,\epsilon )    \}      $, and $\tau' =\{t>\tau| X_t \in \{-3\epsilon,3\epsilon\}             \}$. Now, define a stochastic process $(U_t)_{t \geq \tau}$ started on $\mathscr{F}_{\tau}$ that satisfies ~\eqref{eq:main} with $U_\tau =-2\epsilon$. From Proposition \ref{prop:Sec4Dominance}, we see that $ U_t<X_t, \forall t\in (\tau,\tau')    $. Now, we can see that $ \mathbb{P} (\tau ' =\infty)= 0$. Indeed, $   \mathbb{P} (\tau ' =\infty) \leq \mathbb{P} ( U_t \leq 3\epsilon \forall t\geq \tau   )\leq 1-\mathbb{   P}  (U_t \rightarrow \infty  ) =0 $. 
			$\hfill \blacksquare$			
			
			\textbf{Proof Theorem \ref{thm:genconv} part \ref{enum:conv1}:} Suppose $\mathcal{N}=(-3\epsilon, 3\epsilon)$ for $\epsilon>0$. Without loss of generality and for the purposes of this proof,  assume that $\epsilon<\min \left (\left ( 
			\frac{C(1)}{k}  \right )^{\frac{1}{k-1}}  , \left ( 
			\frac{C(c)}{k}  \right )^{\frac{1}{k-1}}   \right )$. Pick a time $z$ such that $h(t)\geq -\frac{3}{2}\left ( 
			\frac{k}{C(c)}  \right )^{\frac{1}{k-1}}\epsilon$ for all $ t\geq z $. 
			From Proposition \ref{prop:hittingeverypoint}, $\tau<\infty$ with positive probability.
			Now, we define two stochastic processes $(Y_t)_{t\geq \tau},(Y'_t)_{t \geq \tau}$ in the same probability space as $X_t$ started on $\mathscr{F}_\tau$, that satisfy \eqref{eq:mainwithScale} with drift constant $1$ and $c$ respectively. From Proposition \ref{prop:Sec4Dominance}, we see that if $Y_\tau>X_\tau>Y'_\tau$, then $Y_t>X_t>Y'_t$ for all $ t \in (\tau, \tau')$. We set $Y'_\tau$, such that 
			$X_\tau>Y'_\tau     $, and $Z^{Y'}_t=\frac{Y'_t}{h(t)} >\max \left (-\left ( 
			\frac{C(1)}{k}  \right )^{\frac{1}{k-1}}  , -\left ( 
			\frac{C(c)}{k}  \right )^{\frac{1}{k-1}}   \right )$. Now, we should show that $\{\tau'=\infty\}\cap\{Y_t\rightarrow 0\} \cap \{Y'_t\rightarrow 0\} $ is non trivial. Take $\epsilon _1$ and $\epsilon _c$, both less than $\epsilon$, as in the statement of Lemma \ref{lemma:Main3.3} for $Y_t$ and $Y'_t$ respectively, and pick $\epsilon '= \min (\epsilon _1, \epsilon _c   )$. For $\epsilon '$, using Lemma \ref{lemma:Main3.3}, we know we can find $\delta_1$ and $\delta_c$ such that on $A_1= \{G_t \in (-\epsilon', \epsilon'    )\,\mathrm{for~all}~ t \in (s, s+\delta_1), \text{and}\,\, G_t \in (-2\epsilon'   ,-\frac{9}{10}\epsilon'    ) \,\mathrm{for~all}~ t \in  (s+\delta_1,\infty)\}   $ we have $Y_s\rightarrow 0$, and on
			$A_c= \{G_t \in (-\epsilon', \epsilon'    )\,\mathrm{for~all}~ t \in (s, s+\delta_c), \text{and}\,\, G_t \in (-2\epsilon'   ,-\frac{9}{10}\epsilon'    ) \,\mathrm{for~all}~ t \in  (s+\delta_c,\infty)\}   $ we have $Y' _s\rightarrow 0$. From here, since $A\cap A'$ is non trivial, we only need to argue that $ \{\tau'=\infty\}   \supset A\cap A'   $. From the remark of Lemma \ref{lemma:Main3.3} we see that $Y_t$ and $Y'_t$ always stay below 0 on  
			$A\cap A'$. Also, from Proposition \ref{prop:StayingClose}, we see that $Z^{Y'}_t > -2\left ( 
			\frac{C(c)}{k}  \right )^{\frac{1}{k-1}} $. Equivalently, and using that $h(t)\geq -\frac{3}{2}\left ( 
			\frac{k}{C(c)}  \right )^{\frac{1}{k-1}}\epsilon$,
			\begin{align*}
			Y'_t &> 2h(t)\left ( \frac{C(c)}{k}  \right )^{\frac{1}{k-1}} \\
			&\geq -3\epsilon
			\end{align*} 
			$\hfill \blacksquare$ 
			
			We now prove the second part of Theorem \ref{thm:genconv}.
			
			\textbf{Proof Theorem \ref{thm:genconv} part \ref{enum:conv2}:} Let $\mathcal{N}= (-3\epsilon, 0).$ Define $\tau=\inf \{t | X_t \in (-2\epsilon,-\epsilon)       \}$, and the exit time from $\mathcal{N}$, $\tau_e=\inf \{t | X_t \not \in (-3\epsilon,0)       \}$. From Proposition \ref{prop:hittingeverypoint}, we have that $\tau <\infty$ holds with positive probability. Define $(Y_t)_{\tau \leq t\leq \tau_e },\,(Y_t)_{\tau \leq t\leq \tau_e}$ to be two processes that satisfy \eqref{eq:absval} with constants $k_1,k_2$ respectively. Suppose that $Y_\tau<X_\tau<Y'_\tau$ and $Y_\tau,Y'_\tau \in (-2\epsilon,\epsilon)$. Then from Proposition  \ref{prop:Sec4Dominance}, we get $Y_t<X_t<Y'_t, \,\mathrm{for~all}~ t \in (\tau,\tau_e)$. Now, using Proposition \ref{prop:comparisonk=1}, there is an event  $A$ such that $Y_t,Y'_t \in (-3\epsilon,0), \,\mathrm{for~all}~ t\geq \tau$. Consequentially $X_t\in (-3\epsilon,0), \,\mathrm{for~all}~ t\geq \tau$ since $\tau_e=\infty$ on $A$. Finally, using Lemma \ref{lem:ContLInf} we conclude that $Y_t\rightarrow 0$ on $A$, hence also $X_t\rightarrow 0$ on $A$. 
			$\hfill \blacksquare$
			
		\end{section}
		
		\begin{section}{The discrete model}	
			\begin{subsection}{Analysis of $X_t$ when $\frac{1}{2} + \frac{1}{2k} > \gamma $, $k>1$ and $\gamma\in(1/2,1)$     }\label{sbec:DiscNconv}
				Before proving Theorem \ref{thm:gennonconvdisc}, as described in  section \ref{sbec:DiscIntro}, we assume that $X_n$ satisfies 
				\begin{equation}\label{eq:discnongenconv}
				X_{n+1} - X_n\geq\frac{|X_n|^k}{n^{\gamma}} +  \frac{Y_{n+1}}{n^{\gamma}},~k>1~\text{and}~\gamma \in (1/2,1),
				\end{equation}
				where $Y_n$ are a.s. bounded and $E(Y_{n+1}|\mathscr{F}_n)=0$. In this section we additionally require $Y_n$ to satisfy 
				$   E(Y_{n+1}^2|\mathscr{F}_n    )\geq l>0 $.
				\begin{thm}\label{thm:discnonconvmain}
					Let $(X_n)_{n\geq 1}$ solve ~\eqref{eq:discnongenconv}. When  $1/2 + 1/2k >\gamma$, $X_t\rightarrow \infty$ a.s.
				\end{thm}
				Now, we develop the necessary tools to prove this theorem. 
				\begin{prop}\label{propClstoZ}
					The process $(X_n)_{n\geq 1}$ gets close to the origin infinitely often. More specifically, for $\beta<0$ the event $\{X_n \geq \beta n^{\frac{1-2\gamma}{2}   } \text{i.o.}  \} $ has probability $1$.
				\end{prop}
				\textbf{Proof:} Now, from the restrictions on $k$ we obtain
				\begin{align*}
				\frac{1}{2} +\frac{1}{2k} \geq \gamma &\iff \frac{1-k}{2k} \geq \gamma -1 \iff \frac{k-1}{k} \leq 2(1-\gamma )\\
				\frac{k-1}{1-\gamma} \leq \frac{2   }{2\gamma -1} &\iff
				\frac{(1-\gamma )}{k-1}  \geq \frac{2\gamma-1   }{2}.
				\end{align*} 
				Set $h(t)  = -t^{ \frac{1-\gamma}{1-k}     }$, and define 
				$Z_n   = -\frac{X_n}{h_n}       $. From here we get the following recursion,
				\begin{align*}
				Z_{n+1}-Z_{n}&\geq-\frac{X_{n+1}      }{h(n+1)} +\frac{X_{n}      }{h(n)}\\
				&\geq-X_n( \frac{1}{h(n+1)}    - \frac{1}{h(n)}       )  - \frac{ |X_n|^k    }{n^\gamma h(n+1)    } - \frac{ Y_{n+1}  }{n^\gamma h(n+1)} \\
				&\geq X_n \frac{1-\gamma}{ k-1   }  \xi_n ^{ -\frac{1-\gamma}{ 1-k   }-1         }-
				\frac{ |X_n|^k    }{n^{\gamma} h(n+1)    }-\frac{ Y_{n+1}  }{n^\gamma h(n+1)}\\
				&\geq   \frac{   X_n  }{h(n+1) n^\gamma}  \left( \frac{1-\gamma}{ k-1   }  \xi_n ^{ -\frac{1-\gamma}{ 1-k   }-1         }  h(n+1)n^\gamma-\frac{|X_n|^k }{X_n}               \right)  -\frac{ Y_{n+1}  }{n^\gamma h(n+1)} \\
				&\geq \frac{   X_n  }{h(n+1) n^\gamma}  \left (-a_n \frac{1-\gamma}{ k-1   } |h(n)|^{k-1}-\frac{|X_n|^k }{X_n}               \right)  -\frac{ Y_{n+1}  }{n^\gamma h(n+1)}\\
				&\geq\frac{   X_n  }{h(n+1) n^\gamma}  \left (- \frac{2(1-\gamma)}{ k-1   } |h(n)|^{k-1}-\frac{|X_n|^k }{X_n}               \right )  -\frac{ Y_{n+1} }{n^\gamma h(n+1)}\\
				\end{align*}
				where $$a_n =      \frac{   \xi_n ^{\frac{    -(1-\gamma)}{1-k} -1     } h(n+1)n^{\gamma} }{ -|h(n) |^{k-1}          }  .$$ 
		        It is easy to verify that 
				$a_n \rightarrow 1$, whence the last inequality in the previous calculation for large enough $n$. 
				
				Define ${G'}_{s,n}=\sum_{ i=s}^{n-1}  \frac{Y{i+1}}{i^{\gamma}h(i+1)   } $, we will see that ${G'}_{1,n}$ grows big enough so that  $Z_n$ must get, at certain times, close enough to the origin so that $X_n$ surpasses a constant multiple of $h(n)$. To this end, we have the following lemma,
				\begin{lemma}\label{DivDiscNoise}
					$\limsup _{ n \rightarrow \infty}   {G'}_{1,n} = \infty$ a.s. 
				\end{lemma}
				We use the following theorem, for a reference see~\cite{MR1159567} page 676 Theorem 1.,
				\begin{thm}
					Let $X_n$ be a martingale difference such that $E(X_i^2|\mathscr{F}_{i-1}) <\infty$. Set $s_n^2=\sum_{i=1}^{n}E(X_i^2|\mathscr{F}_{i-1})$, and define $\phi(x) = ( 2\log _2(x^2\vee e^2)    )^{\frac{1}{2}   }$. We assume that $s_n\rightarrow \infty$ a.s. and that $|X_i|\leq \frac{K_i s_i}{\phi(s_i)}$ a.s. where $K_i$ is $\mathscr{F}_{i-1}$ measurable with $ \limsup_{i \rightarrow \infty} K_i <K$ for some constant $K.$ Then there is a positive constant $\epsilon(K)$ such that 
					$\limsup_{n \rightarrow \infty}\sum_{i=1}^{n} \frac{X_i}{s_n\phi(s_n) }\geq \epsilon(K) $ a.s.
				\end{thm} 
				It is clear that ${G'_{ 1,n     }}$ satisfies all the hypothesis required for the aforementioned theorem to hold.
				$\blacksquare$
				
				From Lemma \ref{DivDiscNoise}, it is immediate that for any random time $s$ (not necessarily a stopping time) $\limsup_{n\rightarrow \infty} {G'}_{s,n}  = \infty  $, a.s.  
				
				Now, we return to prove proposition \ref{propClstoZ}. Assume that there is $n_0,$ such that $   X_n < -(\frac{3(1-\gamma)}{k-1})^{\frac{1}{k-1}} n^{  \frac{1-\gamma}{1-k}   }    $, for all $n\geq n_0$. Then, since $  \frac{|x|^k}{x}$ is increasing  we get that 
				$   \frac{|X_n|^k}{X_n} <-\frac{3(1-\gamma)}{k-1} n^{   1-\gamma    }   $. 
				Therefore,  
				$$   -\frac{2(1-\gamma)}{ k-1   } |h(n)|^{k-1}-\frac{|X_n|^k }{X_n}           >   -\frac{2(1-\gamma)}{ k-1   } n^{1-\gamma}+ \frac{3(1-\gamma)}{k-1} n^{   1-\gamma    }  = \frac{(1-\gamma)}{k-1} n^{   1-\gamma    }  >0      $$
				
				So,
				\begin{align*}
				Z_n&\geq Z_{n_0} +\sum _{i=n_0}^{n}\frac{   X_i  }{h(i+1)i^\gamma}  ( \frac{2(1-\gamma)}{ k-1   } |h(i)|^{k-1}-\frac{|X_i|^k }{X_i}               )+{G'}_{n_0,n} \\
				&>  Z_{n_0}+   {G'}_{n_0,n},   
				\end{align*} 
				which gives $ \limsup _{n \rightarrow \infty}Z_n = \infty $ which is a contradiction since this would imply $X_n \geq 0$, infinitely often.
				
				Since $ n^{\frac{1-\gamma}{1-k}} = o(    n^{\frac{1-2\gamma}{2}    })$, for every constant $\beta<0$, the event $\{ X_n \geq \beta n^{ \frac{1-2\gamma}{ 2}    }  \text{\,i.o.}   \}  $, holds a.s.
				$\hfill \blacksquare$			
				\begin{lemma}
					For any $n$, we can find $a_1>0,\delta>0$ such that $   \mathbb{P}(\sup_{u\geq n}G_{{n},u}\geq a_1 n^{ \frac{1-2\gamma}{2}}  |\mathscr{F}_n    )> \delta $ and 
					$   \mathbb{P}(G_{n,\infty}\geq a_1 n^{ \frac{1-2\gamma}{2}}  |\mathscr{F}_n    )>\delta$.
				\end{lemma}
				\textbf{Proof:} Define $\tau = \inf \{u\geq n| G_{n,u} \notin (-a_2n^{\frac{1-2\gamma}{2}},a_2n^{\frac{1-2\gamma}{2}}   )      \}    $. We calculate the stopped variance of $ G_{\tau}:=G_{n,\tau}$. We will do so recursively; fix $m\geq n$ and calculate, 
				\begin{align*}
				E (( G_{ \tau \wedge {m+1}      }   )  ^2  |\mathscr{F}_n )-E ( G_{\tau \wedge m}   )^2  |\mathscr{F}_n )&=   E \left (1_{\tau>m} \left ( 2\frac{Y_{m+1}}{m^\gamma    } G_{m}             +\frac{   Y_{m+1}^2}{m^2}     \right )    |\mathscr{F}_n \right )          \\
				&=E \left (1_{\tau>m}  2\frac{Y_{m+1}}{m^\gamma    } G_{m}|\mathscr{F}_n \right )            +E\left (1_{\tau>m}\frac{   Y_{m+1}^2}{m^2}         |\mathscr{F}_n \right )    \\
				&=0 + E\left (1_{[\tau>m]} E\left(\frac{   Y_{m+1}^2}{m^{2\gamma}}     |\mathscr{F}_{m}  \right)  |\mathscr{F}_n   \right )\\
				&\geq \epsilon \frac{1}{m^{2\gamma}} E(1_{[\tau>m]} |\mathscr{F}_n    )\\
				&\geq \epsilon \frac{1}{m^{2\gamma}} \mathbb{P}(\tau=\infty |\mathscr{F}_n )   .
				\end{align*}
				Therefore,
				\begin{align} \label{VarIn}
				E (( G_{ \tau \wedge {m}      }   )  ^2  |\mathscr{F}_n )&\geq E ( (G_{\tau \wedge n }  )^2  |\mathscr{F}_n )+c\mathbb{P}(\tau=\infty |\mathscr{F}_n ) (n^{1-2\gamma}        -(m-1)^{1-2\gamma}   )\nonumber\\
				&= c\mathbb{P}(\tau=\infty |\mathscr{F}_n ) (n^{1-2\gamma}        -(m-1)^{1-2\gamma}   ).
				\end{align}
				Notice that since $Y_n$ are a.s. bounded, $|G_\tau| \leq a_2 n^{\frac{1-2\gamma}{2}   }   +\frac{M}{n^{\gamma}}$, and since $\gamma > \gamma-1/2$,
				we get that $|G_\tau| \leq 2a_2n^{\frac{1-2\gamma}{2}   } $ for $n$ large enough. For $m$ large, we can find a constant $c'$ such that $n^{1-2\gamma}        -(m-1)^{1-2\gamma}  \geq c'n^{1-2\gamma}$. Using \eqref{VarIn}, we obtain
				$$ \frac{ 2a_2n^{\frac{1-\gamma}{2}   } }{\epsilon c'n^{1-2\gamma}}\geq    \mathbb{P}(\tau=\infty |\mathscr{F}_n ). $$
				Choosing $a_2$ small enough we may conclude $\mathbb{P}   (      \tau < \infty | \mathscr{F}_n) > 1/4 $, for all $n$ large enough.
				
				Now, we take any martingale $M_n$ starting at 0, such that it exits an interval $ (-2a,2a)    $, with at least probability $p$, and $|M_{n+1}-M_n|< a $, a.s. Then, we stop the martingale upon exiting the interval $(-2a,2a)$; namely, define $\tau_{-}$ to be the first time $M_n$ goes below $-2a$ and $\tau_{+}$ to be the first time that $M_n$ surpasses  $2a$, and set $\tau = \tau_{-} \wedge \tau_{+}$. Using the optimal stopping theorem for the bounded martingale $M_{\tau \wedge t}$ and taking $t$ to infinity gives
				\begin{align*}
				0=E(M_{\tau})&\leq 	\mathbb{P}(\tau_{-}< \tau_{+}    )(-2a) +\mathbb{P}(\tau_{-}> \tau_{+}    )(3a) +\mathbb{P}(\tau =\infty)2a\\
				&=-2ap +(1-p)2a+\mathbb{P}(\tau_{-}> \tau_{+}    )(5a) .
				\end{align*} 
				So, $   \mathbb{P}(\tau_{-}> \tau_{+}    ) \geq   \dfrac{4p-1}{5}       $, which implies that $\mathbb{P}(\sup M_n \geq 2a ) \geq \dfrac{4p-1}{5}$.
				
				The previous applied to $G_{n,u}$ given $\mathscr{F}_n$, concludes the lemma. Indeed, since the probability $p$, of exiting the interval is bigger than $1/4$, we may deduce that  $\dfrac{4p-1}{5}>0$.
				
				For the second part of the lemma, we use the following inequality: let $M_n$ be a martingale such that $M_0=0$ and $E(M_n^2) <\infty$. Then 
				$\mathbb{P} (  \max_{n\geq u \geq 0} M_u \geq \lambda    ) \leq \frac{  E( M_n^2 )}{E( M_n^2 )+ \lambda^2   }$ (for a reference see \cite{dur13}, page 213, exercise 5.4.5). Let $\tau$ be the first time $G_{n,u}$, surpasses $ a_2 n^{1-2\gamma} $. Condition on $[\tau < \infty]$, and notice that $G_{n,\infty}\geq a_2n^{1-2\gamma}$ when $\inf _{ u\geq \tau}G_{\tau,u} > -\frac{a_2}{2}n^{1-2\gamma}$. Using the previous inequality, and the fact that $\frac{x}{x+1}$ is increasing gives
				\begin{align*}
				\mathbb{P}(G_{n,\infty}\leq   \frac{a_2}{2}n^{\frac{1-2\gamma}{2}}  |\mathscr{F}_{\tau},[\tau <\infty]    ) &\leq     \mathbb{P}( \inf _{ u\geq \tau}G_{\tau,u} \leq -\frac{a_2}{2}n^{1-2\gamma}  |\mathscr{F}_{\tau},[\tau <\infty]    ) \\     
				&\leq      \frac{  E((G_{\tau,\infty} )^2 |\mathscr{F}_{\tau},[\tau <\infty] )}{E((G_{\tau,\infty} )^2 |\mathscr{F}_{\tau},[\tau <\infty])+  \frac{a_2^2}{4}n^{1-2\gamma     }}      \\
				&\leq\frac{ c\tau ^{1-2\gamma    }            }{ c\tau ^{1-2\gamma    }  + \frac{a_2^2}{4}n^{1-2\gamma     }    } \\
				&\leq  \frac{ c            }{ c  + \frac{a_2^2}{4}  }. 
				\end{align*}
				Therefore, $\mathbb{P}(G_{n,\infty}\geq \frac{a_2}{2} n^{ \frac{1-2\gamma}{2}}  |\mathscr{F}_n    ) \geq \mathbb{P} (\tau <\infty) \frac{  \frac{a_2^2}{4}  }      { c  + \frac{a_2^2}{4}  }  $, which concludes the lemma.
				$\hfill \blacksquare$
				
				So, for any stopping time $\sigma$, we get the following version of the previous lemma:
				\begin{lemma}
					For any $n$, we can find $a_1>0,\delta_1>0,\delta_2>0$ such that $   \mathbb{P}(\sup_{u\geq \sigma}G_{\sigma,u}\geq a_1 \sigma^{ \frac{1-2\gamma}{2}}  |\mathscr{F}_\sigma    )> \delta_1 $ and 
					$   \mathbb{P}(G_{\sigma,\infty}\geq a_1 \sigma^{ \frac{1-2\gamma}{2}}  |\mathscr{F}_\sigma    )>\delta_2$.
				\end{lemma}
				$\hfill \blacksquare$
				\begin{cor}
					The event $\{   X_n \geq 0\, \text{i.o.}\}$ holds a.s.
				\end{cor}
				\textbf{Proof:} For any $m,n$ we get the lower bound 
				$ X_m-X_n \geq G_{n,m} $. Now, we define an increasing sequence of stopping times $\tau _n$, going to infinity a.s., such that $X_{\tau_n}    \geq \beta \tau_n^{\frac{1-2\gamma}{2}   }$ for $|\beta|<a_1$. From Proposition \ref{propClstoZ}, we can do so, with all $\tau_n$ a.s. finite. Hence, $\mathbb{P}(\sup_{\infty \geq u\geq {\tau_n}} X_m-X_{\tau_n}\geq a_1{\tau_n}^{  \frac{1-2\gamma}{2}}|\mathscr{F}_{\tau_n}  ) \geq \mathbb{P}(\sup_{\infty \geq u\geq {\tau_n} } G_{\tau_n,u } \geq a_1 {\tau_n}^{  \frac{1-2\gamma}{2}    }|\mathscr{F}_{\tau_n}  )>\delta_1 >0 $. Therefore, by Borel-Cantelli on the event $\{X_{\tau_n} \geq  \beta \tau_n^{\frac{1-2\gamma}{2}  }  \}$, we get $\{X_{\tau_n} \geq   0 ~\text{i.o.}  \}$. Therefore $\{X_n \geq  0~\text{i.o.}  \}$ holds a.s.
				$\hfill \blacksquare$
				
				\textbf{Proof of Theorem \ref{thm:discnonconvmain}:} Define $\tau_n$, as in the proof of the previous corollary, such that $X_{\tau_n} \geq 0$. 
				Since $   \mathbb{P}(G_{\tau_n,\infty}\geq a_1 \tau_n^{ \frac{1-2\gamma}{2}}  |\mathscr{F}_{\tau_n}   )>\delta_2$, an application of Borel-Cantelli shows that $\{X_n\geq \frac{a_1}{2}n^{ \frac{1-2\gamma}{2}}  ~\text{i.o.}   \}$ holds a.s. We claim a.s. there are constants $c(\omega)>0$ $m(\omega) $ such that $\{ X_n > c \mathrm{~for~all~} n\geq m  \}=\{\liminf_{\rightarrow \infty} X_n> 0\}$. Indeed, if we define $\tau_0=0$ and $\tau_{n+1} = \{m>\tau_n+1| X_m \geq \frac{a_1}{2} m^{ \frac{1-2\gamma}{2}}  \}$ we see that $\tau_n<\infty$ a.s. and $\tau_n \rightarrow \infty.$ This gives a corresponding filtration, namely $\mathscr{F}_n = \sigma( \tau_n)$.
				
				To finish the claim, we show that $A=\{\liminf_{\rightarrow \infty} X_n\leq 0\}$ has probability zero. To do so, it is sufficient to argue that there is a $\delta$ such that $ \mathbb{P}( A|\mathscr{F}_n   )<1-\delta $ a.s. for all $n\geq 1$. This is immediate from the previous calculation. Indeed,
				\begin{align*}
				\mathbb{P}( A|\mathscr{F}_n   ) & \leq  1-\mathbb{P}\left ( \liminf _{n } X_n \geq 
				\frac{3a_1}{2}\tau_n^{ \frac{1-2\gamma}{2}}  |\mathscr{F}_n \right )\\
				&= 1-\mathbb{P} \left ( \liminf _{n } X_n-\frac{a_1}{2} \tau_n^{ \frac{1-2\gamma}{2}} \geq
			a_1\tau_n^{ \frac{1-2\gamma}{2}}  |\mathscr{F}_n \right ) \\
			&\leq 1-\mathbb{P} ( \liminf _{n }G_{\tau_n,n}  \geq
			a_1\tau_n^{ \frac{1-2\gamma}{2}}  |\mathscr{F}_n )\\
				&<1-\delta_2.
				\end{align*}				
				The process $G_{m,\infty}$ is a.s. finite, and since the drift term $\sum_{ i\geq n } \frac{|X_i|^k}{i^{\gamma}} \rightarrow \infty$, we get that $X_n \rightarrow \infty $.
				$\hfill \blacksquare$

				\textbf{Proof of Theorem \ref{thm:gennonconv}:} We define, $\tau= \{n| X_n \in (-\epsilon, \epsilon)\}$, and $\tau ' = \{\tau'>\tau | X_n \not \in (-\epsilon ' , \epsilon ')   \}$. When $\epsilon$ is small enough, we may assume that $\tau<\infty$ with positive probability, otherwise we have nothing to prove. On $\{\tau<\infty   \}$, couple $X_n$ with $X'_n$, so that $\mathbb{P}(X_n=X'_n,~\tau \leq n\leq \tau ' |[\tau <\infty] )=1  $, where $X_n'$ is a process that solves \eqref{eq:discnongenconv}. Since $X'_n\rightarrow \infty$, a.s., we have that $\tau' <\infty$ a.s. Thus, on $\{\lim_{n \rightarrow \infty} X_n  =0\} $ by Borel-Cantelli implies $\{ X_n=0~\text{i.o.} \}$. Therefore, $\mathbb{P}(\lim_{n \rightarrow \infty} X_n  =0) = 0$.
				$\hfill \blacksquare$
			\end{subsection}

			\begin{subsection}{Analysis of $X_t$ when $\frac{1}{2} + \frac{1}{2k} < \gamma $, $k>1$ and $\gamma\in(1/2,1)$  }\label{sbec:disConv}
				Before proving the main Theorem \ref{thm:genconvdisc2}, as described in the section \ref{sbec:DiscIntro} we will study a process $(X_n)_{n \geq 1}$ that satisfies
				\begin{equation}\label{eq:discconvnongen}
				X_{n+1} - X_n\leq\frac{f(x)}{n^{\gamma}} +  \frac{Y_{n+1}}{n^{\gamma}},\,\gamma \in (1/2,1),\,k\in (1,\infty),
				\end{equation}
				where, $f(x) \leq |x|^k$ when $  x\in(-\epsilon,\epsilon)   $, and $f(x) = |x|^k$ when $  x\in\mathbb{R}\setminus (-\epsilon,\epsilon)   $. Let $x_0<0$, such that $f(x)>M,~\forall x\leq x_0$. We will use $x_0$ in the next lemma.
				
				\begin{lemma}
					Take $C= \max( M,|X_1|,|x_0|     )$. Then $X_n> -2C$ for all $n$, a.s.
				\end{lemma}
				\textbf{Proof:} We can show this by induction. Of course $X_1> -2C$. For the inductive step, we distinguish two cases. First, assume that $-2C<X_n<-C$. Then 
				\begin{align*}
				X_{n+1} &= X_n +\frac{f(X_n)}{n^\gamma} + \frac{Y_{n+1}}{n^{\gamma}} \\
				&\geq   -2C +  \frac{f(X_n)}{n^\gamma} -\frac{M}{n^{\gamma}}\\
				&> -2C. 
				\end{align*}
				Now, assume $X_n\geq -C$. Then
				\begin{align*}
				X_{n+1} &= X_n +\frac{|X_n|^k}{n^\gamma} + \frac{Y_n}{n^{\gamma}} \\
				&\geq  -C + 0 -\frac{M}{n^{\gamma}}\\
				&> -2C.
				\end{align*}
				$\hfill \blacksquare$
				
				 Pick $\epsilon >0$ such that $ \epsilon \leq \min   (  \frac{1}{4}, \frac{1}{2} (\frac{1-\gamma}{3(k-1)})^{ \frac{1}{k-1}   }    )    $. Let $a_n$, be defined as in the previous section.  
				
				\textbf{Claim:} \textit{We can find $n_0$, that satisfies the following properties }
					\begin{enumerate}
						\item $a_n >1/2,\, n\geq n_0$ a.s.
						\item if   $-\frac{X_{n+1}}{h(n+1)}>-2\epsilon$, and $-\frac{X_{n}}{h(n)}\leq-2\epsilon$, then
						$-\frac{X_{n+1}}{h(n+1)} <-\epsilon$, when $n\geq n_0$.
						\item $\mathbb{P} ({G'}_{n_0,n} \in (\frac{-\epsilon}{2}, \frac{\epsilon}{2})\forall n \geq n_0 |\mathscr{F}_{n_0}  ) >0$.
					\end{enumerate} 
				
					\textbf{Proof:} 
						\begin{enumerate}
					\item This is is trivial.
					
					\item  Since $|Y_n|<M$, and $X_n>C$ a.s., then whenever $X_n<0$, we have $|X_{n+1}-X_{n}|=O(n^{-\gamma})$. Also, $ n^{-\gamma} = o(h(n)   )$, since $\gamma > \frac{1-\gamma}{k-1}$. Indeed, $\gamma > \frac{1-\gamma}{k-1}$ is equivalent to $    \gamma >1/k = 1/2k+1/2k$, however $ 1/2>1/2k$ and since $\gamma > 1/2 + 1/2k$ we conclude. Furthermore, notice that $\dfrac{h(n)}{h(n+1)}\rightarrow 1$. 
					
					Calculate 
					\begin{align*}
					-\frac{X_{n+1}}{h(n+1)} &= -\frac{X_{n+1}-X_{n}    }{h(n+1)} -\frac{X_{n}}{h(n)}\cdot\frac{h(n)}{h(n+1)}\\
					&\geq o(1)-2\epsilon \frac{h(n)}{h(n+1)}
					\end{align*}
					Since the $o(1)$ term and  $\frac{h(n)}{h(n+1)}$ depend only on $n$, we conclude 2.
					
					\item Using Doob's inequality, and the fact that $  m^{\gamma} h(m+1)\sim m^{\frac{1-\gamma}{k-1} - \gamma  }   \leq m^{\frac{1-\gamma}{k-1} - \gamma  }\leq m^{\frac{-1-\delta}{2}} $ for some $\delta>0$, we have:
					\begin{align*}
					\mathbb{P}\left ( \sup _{u\geq n_0     } ({G'}_u^{n_0} |\mathscr{F}_{n_0})^2 \geq \frac{\epsilon^2}{4}\right ) &\leq \sum _{m\geq n_0   } \frac{E(Y_{m+1}^2|\mathscr{F}_{n_0} )    }{m^\gamma h(m+1)    }\\
					&\leq     C\sum _{m\geq n_0   } \frac{1    }{m^\gamma h(m+1)    }\\
					&=\sum _{m\geq n_0   } \Theta (m^{ \frac{1-\gamma}{k-1} - \gamma       }    )\\
					&=\sum _{m\geq n_0   } \Theta (m^{-1-\delta}    )\\
					&=  \Theta ({n_0}^{-\delta}    )\rightarrow 0.
					\end{align*}
					\end{enumerate}
					$\hfill \blacksquare$
					
					Notice, that the previous claim holds for any stopping time $
					\tau$, in place of $n$. So, we obtain a version of the previous lemma for stopping times.
					\begin{lemma}\label{lemma:disclittlenoise}
					Let $\tau$ be a stopping time such that $\tau \geq n_0$, where $n_0$ is the same as in the previous claim. Then,
					$\mathbb{P} ({G'}_{\tau,n} \in (\frac{-\epsilon}{2}, \frac{\epsilon}{2})\forall n \geq \tau | \mathscr{F}_{\tau}   ) >0$
					\end{lemma}
					$\hfill \blacksquare$
					
					Let $ \epsilon \leq \min   (  \frac{1}{4}, \frac{1}{2} (\frac{1-\gamma}{3(k-1)})^{ \frac{1}{k-1}   }    )    $, and define a stopping time $\tau=\inf \{n\geq n_0 |  Z_n < -2\epsilon          \}    $.
					\begin{prop}\label{prop:centraldiscconv}
					Let $(X_n)_{n\geq 1}$ that satisfies \eqref{eq:discconvnongen}. When $\tau<\infty$, with positive probability, then $ \mathbb{P}(X_n\rightarrow 0 )>0$. More specifically, the process $(X_n : n\geq \tau)$ converges to zero with positive probability.
					\end{prop}
					\textbf{Proof:}
					We use the expression for $Z_n = -\dfrac{X_n}{h(n)}$, 
					\begin{align*} Z_{n+1}-Z_n  &\leq \frac{   X_n  }{h(n+1) n^\gamma}  (-a_n \frac{1-\gamma}{ k-1   } |h(n)|^{k-1}-\frac{|X_n|^k }{X_n}               )  -\frac{ Y_{n+1}  }{n^\gamma h(n+1)}\\
					&< \frac{   X_n  }{h(n+1) n^\gamma}  (- \frac{1-\gamma}{ 2(k-1)   } |h(n)|^{k-1}-\frac{|X_n|^k }{X_n}               )  -\frac{ Y_{n+1}  }{n^\gamma h(n+1)}.
					\end{align*}
					Set $D_n=\dfrac{   X_n  }{h(n+1) n^\gamma}  (- \dfrac{1-\gamma}{ 2(k-1)   } |h(n)|^{k-1}-\dfrac{|X_n|^k }{X_n}               ) $. Then we have
					\begin{equation} \label{in:DiscMain} 
					Z_m - Z_{\tau}\leq   \sum _{ i = \tau   }^{m-1} D_i + {G'}_{\tau,m},
					\end{equation}
					which we obtained in the previous subsection. 
					
					Now, we will show, by contradiction, that on the event $A=\{{G'}_{\tau,n} \in (\frac{-\epsilon}{2}, \frac{\epsilon}{2}),\forall n \geq \tau \}$ the process satisfies $X_n < 0$ for all $n \geq \tau$. 
					Define $\tau_0 =\inf \{n\geq \tau|  Z_n \geq 0   \}$, and $\sigma = \sup \{
					\tau \leq n <\tau_0   | Z_{n-1}\leq -2\epsilon, \,Z_n>-2\epsilon      \} $. Also, when 
					$Z_n\geq -2\epsilon$ we have $X_n \geq 2\epsilon h(n)= -2\epsilon n^{   \frac{1-\gamma}{1-k}} $. So $  \frac{|X_n|^k}{ X_n    } \geq -(2\epsilon)^{k-1} n^{1-\gamma     }    $. Therefore, by the definition of $\epsilon$, we get
					$$- \dfrac{1-\gamma}{ 2(k-1)   } |h(n)|^{k-1}-\dfrac{|X_n|^k }{X_n}   < \left (- \dfrac{1-\gamma}{ 2(k-1)   }+  \dfrac{1-\gamma}{ 3(k-1)   }\right )n^{-1+\gamma}=-\dfrac{1-\gamma}{ 6(k-1)   }   n^{-1+\gamma}        <0. $$
					Hence $D_n <0$ whenever $Z_n\geq -2\epsilon$.	 
					If  $\{\tau_0<\infty \}\cap A$ has positive probability, then $\{\sigma<\infty \}\cap A$ does also. Thus, on $\{\tau_0<\infty \}\cap A$,
					\begin{align*}
					0\leq Z_{\tau_0} &=   Z_{\tau} +   \sum _{ i = \tau   }^{\tau_0-1} D_i + {G'}_{\tau,\tau_0}  \\
					&= Z_{\tau} -Z_{\sigma} +Z_{\sigma} +  \sum _{ i = \tau   }^{\tau_0-1} D_i + {G'}_{\tau,\tau_0}  \\
					&=Z_{\sigma}-{G'}_{\tau,\sigma}  +{G'}_{\tau,\tau_0}+\sum _{ i = \sigma   }^{\tau_0-1} D_i\\
					&<-\epsilon +\frac{\epsilon}{2} +\frac{\epsilon}{2} +0=0
 					\end{align*}
					which is a contradiction.
					
					Now, we can complete the proof of the proposition. On the event $A$, $X_n<0$ for all $n>\tau$, therefore $ \limsup_{ n \rightarrow \infty} X_n \leq 0$ on $A$. However, by Lemma \ref{discInf} we have $\limsup_{ n \rightarrow \infty} X_n \geq 0$ a.s. Therefore, on $A$, $X_n\rightarrow 0$. $\hfill \blacksquare$
				
					Remark: On $A$ we showed that $X_n$ converges to zero, since for all $n\geq \tau$, $X_n<0$ and the only place to converge is the origin.					
								
					\textbf{Proof of Theorem \ref{thm:genconvdisc2}:} We define $\tau= \{n\geq n_0| X_n      \in (-\epsilon _2,-\epsilon_1    )                 \}$, where $n_0$ is the same as in Lemma \ref{lemma:disclittlenoise}, and $\tau_e=\inf\{n|  X_n\not \in (-3\epsilon, 3\epsilon)          \}$. 
					Let $(X'_n: n\geq \tau)$ be a process that satisfies \eqref{eq:discnongenconv}. Then we couple $(X_n)$ with $(X'_n)$ on $\{ \tau <\infty     \}$ such that 
					$\mathbb{P}(X_n = X'_n, \tau \leq n \leq \tau_e | \{ \tau <\infty     \}      ) =1$. To show that $X'_n$, converges to zero with positive probability, first we need to verify that the conditions for Proposition \ref{prop:centraldiscconv} are met. The only thing we need to check is that $Z'_\tau= - \frac{X'_\tau}{h(\tau)} <-2\epsilon$. However, since $h(t)\rightarrow 0$ this is always possible by choosing $n_0$ large enough. Furthermore, by Proposition \ref{prop:centraldiscconv}, we see that there is an event of positive probability such that $X'_n\rightarrow 0$, where $\tau_e$ is infinite conditioned on this event. Therefore, $X_n$ converges to 0 with positive probability.
					$\hfill \blacksquare$
				\end{subsection}
			\end{section}
			
			\section*{Acknowledgements}
The author would like to thank Marcus Michelen, Albert Chen and Josh Rosenberg for constructive criticism of the manuscript.			
			
			\bibliographystyle{alpha}
			\bibliography{Bib}

		\end{document}